%& amstex
\input amstex
% This is an AMS-TeX file and should be compiled
% using AMS-TeX. 
\magnification=1200
\loadmsbm
\loadeufm
\loadeusm
\UseAMSsymbols
\input amssym.def

\font\BIGtitle=cmr10 scaled\magstep3
\font\bigtitle=cmr10 scaled\magstep1
\font\boldsectionfont=cmb10 scaled\magstep1
\font\section=cmsy10 scaled\magstep1

\def\scr#1{{\fam\eusmfam\relax#1}}
\def\scrA{{\scr A}}
\def\scrB{{\scr B}}
\def\scrC{{\scr C}}
\def\scrD{{\scr D}}
\def\scrE{{\scr E}}
\def\scrF{{\scr F}}
\def\scrG{{\scr G}}
\def\scrH{{\scr H}}
\def\scrI{{\scr I}}
\def\scrL{{\scr L}}
\def\scrK{{\scr K}}
\def\scrJ{{\scr J}}
\def\scrM{{\scr M}}
\def\scrN{{\scr N}}
\def\scrO{{\scr O}}
\def\scrP{{\scr P}}
\def\scrQ{{\scr Q}}

\def\scrS{{\scr S}}
\def\scrU{{\scr U}}
\def\scrR{{\scr R}}
\def\scrT{{\scr T}}
\def\scrV{{\scr V}}
\def\scrX{{\scr X}}
\def\scrY{{\scr Y}}
\def\scrZ{{\scr Z}}
\def\scrW{{\scr W}}
\def\gr#1{{\fam\eufmfam\relax#1}}

	\def\gro{{\gr o}}
	 
\def\grQ{{\gr Q}}	
	
\def\grS{{\gr S}}	\def\grs{{\gr s}}

\def\db#1{{\fam\msbfam\relax#1}}

\def\dbA{{\db A}} 
\def\dbC{{\db C}} 
 \def\dbF{{\db F}}
\def\dbG{{\db G}} \def\dbH{{\db H}}

 \def\dbN{{\db N}}
 
\def\dbQ{{\db Q}} \def\dbR{{\db R}}

 \def\dbZ{{\db Z}}

\def\Ker{\text{Ker}}
\def\der{\text{der}}
\def\Sh{\hbox{\rm Sh}}

\def\s{\text{s}}

\def\m{\text{m}}

\def\Res{\text{Res}}
\def\ab{\text{ab}}

\def\ad{\text{ad}}

\def\Gal{\text{Gal}}
\def\Hom{\text{Hom}}
\def\End{\text{End}}
\def\Spec{\text{Spec}\;}

\def\Lie{\text{Lie}}

\def\leaderfill{\leaders\hbox to 1em
     {\hss.\hss}\hfill}
\def\nspace{\lineskip=1pt\baselineskip=12pt\lineskiplimit=0pt}

     %the way to use this is "\Proclaim{Theorem 1.1.}" for instance.
\def\endproof{$\hfill \square$}
\def\finishproclaim{\par\rm
     \ifdim\lastskip<\medskipamount\removelastskip
     \penalty55\medskip\fi}
\def\proof{\par\noindent {\it Proof:}\enspace}
\def\references#1{\par
  \centerline{\boldsectionfont References}\medskip
     \parindent=#1pt\nspace}
\def\Ref[#1]{\par\hang\indent\llap{\hbox to\parindent
     {[#1]\hfil\enspace}}\ignorespaces}
\def\Item#1{\par\smallskip\hang\indent\llap{\hbox to\parindent
     {#1\hfill$\,\,$}}\ignorespaces}
\def\ItemItem#1{\par\indent\hangindent2\parindent
     \hbox to \parindent{#1\hfill\enspace}\ignorespaces}

\def\Ge{{\mathchoice{\,{\scriptstyle\ge}\,}
  {\,{\scriptstyle\ge}\,}
  {\,{\scriptscriptstyle\ge}\,}{\,{\scriptscriptstyle\ge}\,}}}

\def\arrowsim{\,\smash{\mathop{\to}\limits^{\lower1.5pt
  \hbox{$\scriptstyle\sim$}}}\,}

\def\doublemaprights#1#2#3#4{\raise3pt\hbox{$\mathop{\,\,\hbox to     
#1pt{\rightarrowfill}\kern-30pt\lower3.95pt\hbox to
     #2pt{\rightarrowfill}\,\,}\limits_{#3}^{#4}$}}

\def\rightcapdownarrow{\raise9pt\hbox{$\ssize\cap$}\kern-7.75pt
     \Big\downarrow}

\def\rcapmapdown#1{\rightcapdownarrow\kern-1.0pt\vcenter{
     \hbox{$\scriptstyle#1$}}}

\def\rmapdown#1{\Big\downarrow\kern-1.0pt\vcenter{
     \hbox{$\scriptstyle#1$}}}
\def\rightsubsetarrow#1{{\ssize\subset}\kern-4.5pt\lower2.85pt
     \hbox to #1pt{\rightarrowfill}}
\def\longtwoheadedrightarrow#1{\raise2.2pt\hbox to #1pt{\hrulefill}
     \!\!\!\twoheadrightarrow}

\def\Gal{\operatorname{\hbox{Gal}}}
\def\Hom{\operatorname{\hbox{Hom}}}

\def\im{\hbox{Im}}

\NoBlackBoxes
\parindent=25pt
\document
\footline={\hfil}

\null
\noindent 
\centerline{\BIGtitle Good Reductions of Shimura Varieties of Hodge Type}
\bigskip
\centerline{\BIGtitle  in Arbitrary Unramified Mixed Characteristic, Part II}
\vskip 0.4 cm
\centerline{\bigtitle Adrian Vasiu, Binghamton University}
\vskip 0.4 cm
\centerline{July 24, 2012}

\footline={\hfill}
\vskip 0.5 cm
\noindent
{\bf ABSTRACT.} We prove a conjecture of Milne pertaining to the existence of integral canonical models of Shimura varieties of abelian type in arbitrary unramified mixed characteristic $(0,p)$. As an application we prove for $p=2$ a motivic conjecture of Milne pertaining to integral canonical models of Shimura varieties of Hodge type. 

\footline={\hss\tenrm \folio\hss}
\pageno=1

\bigskip
\noindent
{\boldsectionfont 1. Introduction}
\bigskip

Let $p\in\dbN$ be an arbitrary prime. Let $\dbZ_{(p)}$ be the localization of $\dbZ$ at its prime ideal $(p)$. In this paper we prove the existence of {\it integral canonical models} of {\it Shimura varieties} of {\it abelian type} in unramified mixed characteristic $(0,p)$ (i.e., over finite, \'etale $\dbZ_{(p)}$-algebras). In this introduction we first begin by recalling basic things on {\it Siegel modular varieties} and {\it Mumford moduli schemes}. Then we recall basic types of Shimura varieties and previous works on the existence of integral canonical models. We will end up the introduction by stating our main results and by outlining the strategies to prove them. Let $d\in\dbN^{\ast}$.

\bigskip\noindent
{\bf 1.1. Siegel modular varieties and Mumford moduli schemes.} Let $(A,\lambda_A)$ be a principally polarized abelian variety over $\dbC$ of dimension $d$. Let $L:=H_1(A(\dbC),\dbZ)$ be the first homology group of the analytic space $A(\dbC)$ with coefficients in $\dbZ$; it is a free abelian group of rank $2d$. Let $\psi:L\times L\to \dbZ$ be the perfect alternating form on $L$ induced naturally by $\lambda_A$. Let $W:=L\otimes_{\dbZ} \dbQ=H_1(A(\dbC),\dbQ)$. The classical {\it Hodge theory} provides us with a Hodge decomposition
$$L\otimes_{\dbZ} \dbC=W\otimes_{\dbQ} \dbC=F^{-1,0}\oplus F^{0,-1}\leqno (1)$$
such that under the standard complex conjugation of $W\otimes_{\dbQ} \dbC$ we have an identity $\overline{F^{-1,0}}=F^{0,-1}$. More precisely, one can identify $F^{-1,0}=\Lie(A)$ and $F^{0,-1}=\Hom_{\dbC}(H^1(A,\scrO_A),\dbC)$, where $\scrO_A$ is the structured ring sheaf on $A$. Both $F^{-1,0}$ and $F^{0,-1}$ are isotropic with respect to $\psi$ and in fact $2\pi i\psi$ is a polarization of the Hodge $\dbZ$-structure on $L$ defined by (1). Thus to (1) corresponds naturally a monomorphism 
$$x_A:\dbC^{\ast}\hookrightarrow \pmb{\text{GSp}}(W\otimes_{\dbQ} \dbR,\psi)$$
of reductive groups over $\dbR$. Let $\scrS$ be the $\pmb{\text{GSp}}(W,\psi)(\dbR)$-conjugacy class of $x_A$. If $C_A$ is the centralizer of $x_A$ in $\pmb{\text{GSp}}(W\otimes_{\dbQ} \dbR,\psi)$, then we have $\scrS=\pmb{\text{GSp}}(W\otimes_{\dbR} \dbR,\psi)(\dbR)/C_A(\dbR)$ and thus $\scrS$ gets a natural structure of a hermitian symmetric domain isomorphic to two copies of the Siegel space of genus $d$. The pair $(\pmb{\text{GSp}}(W,\psi),\scrS)$ is called a {\it Siegel modular pair} and it is the most studied type of {\it Shimura pairs}. We have a canonical identification
$$A(\dbC)=F^{0,-1}\backslash (W\otimes_{\dbQ} \dbR)/L\leqno (2)$$
of complex Lie groups. If $(B,\lambda_B)$ is another principally polarized abelian variety over $\dbC$ of dimension $d$, then there exists an element $g\in \pmb{\text{GSp}}(W,\psi)(\dbR)$ such that the complex Lie group $B(\dbC)$ is isomorphic to $g(F^{0,-1})\backslash (W\otimes_{\dbQ} \dbR)/L$ in such a way that the perfect alternating form on $L$ defined by $\lambda_B$ is either $\psi$ or $-\psi$. One gets that the course moduli space of principally polarized abelian varieties over $\dbC$ of dimension $d$ is isomorphic to 
$$\pmb{\text{GSp}}(L,\psi)(\dbZ)\backslash \scrS\leqno (3)$$
(see [BB, Thm. 10.11] for the canonical structure of (3) as a normal, quasi-projective variety over $\dbC$). If $\dbA_f:=\widehat{\dbZ}\otimes_{\dbZ}\dbQ$ is the ring of finite ad\`eles of $\dbQ$, then (3) is isomorphic to 
$$\pmb{\text{GSp}}(W,\psi)(\dbQ)\backslash [\scrS\times (\pmb{\text{GSp}}(W,\psi)(\dbA_f)/\pmb{\text{GSp}}(L,\psi)(\widehat{\dbZ}))].\leqno (4)$$
\indent
Let $N\in \dbN\setminus (p\dbN\cup\{1,2\})$. Let 
$$K(N):=\{g\in \pmb{\text{GSp}}(L,\psi)(\widehat\dbZ)|g\;\text{modulo}\; N\widehat{\dbZ}\;\text{is}\;\text{identity}\}\;\text{and}\;K_p:=\pmb{\text{GSp}}(L,\psi)(\dbZ_p).$$ 
Let $\scrA_{d,1,N}$ be the {\it Mumford moduli scheme} over $\dbZ[{1\over N}]$ that parameterizes isomorphism classes of principally polarized abelian schemes over $\dbZ[{1\over N}]$-schemes that are of relative dimension $d$ and that are endowed with a symplectic similitude level-$N$ structure (cf. [MFK, Thms. 7.9 and 7.10] applied to symplectic similitude level structures instead of simply level structures). The $\dbZ[{1\over N}]$-scheme $\scrA_{d,1,N}$ is smooth and quasi-projective, cf. loc. cit. Similarly to (4) one gets a natural identification 
$$\scrA_{d,1,N}(\dbC)=\pmb{\text{GSp}}(W,\psi)(\dbQ)\backslash [\scrS\times (\pmb{\text{GSp}}(W,\psi)(\dbA_f)/K(N))].\leqno (5)$$
Based on (5) and on classical works of Shimura, Taniyama, etc., one gets an identification
$$\scrA_{d,1,N,\dbQ}=\Sh(\pmb{\text{GSp}}(W,\psi),\scrS)/K(N)\leqno (6)$$
of $\dbQ$--schemes, where $\Sh(\pmb{\text{GSp}}(W,\psi),\scrS)$ is the {\it canonical model} as defined in [De1] of the complex Shimura variety 
$$\Sh(\pmb{\text{GSp}}(W,\psi),\scrS)_{\dbC}=\text{proj.}\text{lim.}_{K\in \Sigma(\pmb{\text{GSp}}(W,\psi))} \pmb{\text{GSp}}(W,\psi)(\dbQ)\backslash [\scrS\times (\pmb{\text{GSp}}(W,\psi)(\dbA_f)/K)].$$
Here $\Sigma(\pmb{\text{GSp}}(W,\psi))$ is the set of compact, open subgroups of $G(\dbA_f)$ endowed with the inclusion relation. Thus 
$$\scrM:=\text{proj.}\text{lim.}_{N\in \dbN\setminus (p\dbN\cup\{1,2\})} \scrA_{d,1,N}$$
is a $\dbZ_{(p)}$-scheme which is an integral model of $\pmb{\text{GSp}}(W,\psi)(\dbQ)\backslash [\scrS\times (\pmb{\text{GSp}}(W,\psi)(\dbA_f)/K_p)]$ over $\dbZ_{(p)}$. In fact $\scrM$ is the integral canonical model of $\pmb{\text{GSp}}(W,\psi)(\dbQ)\backslash [\scrS\times (\pmb{\text{GSp}}(W,\psi)(\dbA_f)/K_p)]$ over $\dbZ_{(p)}$ in the strongest sense of [Va1, Def. 3.2.3 6)] (see [Va1, Ex. 3.2.9] or [Mi2, Thm. 2.10]) and thus also in the weaker sense of [Moo]. From this and [VZ, Cors. 5 and 30] one gets that the regular, formally smooth $\dbZ_{(p)}$-scheme $\scrM$ is uniquely determined by its generic fibre $\scrM_{\dbQ}$ and by the following universal property:

\medskip\noindent
{\it If $Z$ is a regular, formally smooth scheme over $\dbZ_{(p)}$, then each morphism $Z_{\dbQ}\to \scrM_{\dbQ}$ of $\dbQ$--schemes extends uniquely to a morphism $Z\to \scrM$ of $\dbZ_{(p)}$-schemes.}

\medskip
The goal of this paper is to generalize the existence and the uniqueness of $\scrM$ to the case of all Shimura varieties of abelian type (i.e., to all Shimura varieties that are moduli spaces of polarized abelian motives endowed with level structures and motivic tensors). 

\bigskip\noindent
{\bf 1.2. Types of Shimura pairs.} Let $G$ be a reductive subgroup of $\pmb{\text{GSp}}(W,\psi)$ for which any one of the following two equivalent statements holds:

\medskip
{\bf (i)} no simple compact factor of the adjoint group $G^{\ad}$ of $G$ becomes compact over $\dbR$ and there exists an element $x\in\scrS$ which factors through $G_{\dbR}$;

\smallskip
{\bf (ii)} there exists an element $x\in\scrS$ which factors through $G_{\dbR}$ and $G^{\ad}_{\dbR}$ is generated by the $G(\dbR)$-conjugates of the homomorphism $x^{\ad}:\dbC^{\ast}\to G^{\ad}_{\dbR}$ induced by $x:\dbC^{\ast}\to G_{\dbR}$.

\medskip
Let $\scrX$ be the $G(\dbR)$-conjugacy class of an element $x\in\scrS$ that factors through $G_{\dbR}$. The pair $(G,\scrX)$ is a Shimura pair in the sense of [De2] and we have an injective map $f:(G,\scrX)\hookrightarrow (\pmb{\text{GSp}}(W,\psi),\scrS)$ of Shimura pairs. A Shimura pair $(G_1,\scrX_1)$ is called of {\it Hodge type} if it is isomorphic to a Shimura pair of the form $(G,\scrX)$ (for some $d\in\dbN^{\ast}$). Let $\scrX_1^{\ad}$ be the $G_1^{\ad}(\dbR)$-conjugacy class of the homomorphism $x_1^{\ad}:\dbC^{\ast}\to G_{1,\dbR}^{\ad}$ induced by some element $x_1\in\scrX_1$. The pair $(G_1^{\ad},\scrX_1^{\ad})$ is a Shimura pair called the {\it adjoint Shimura pair} of $(G_1,\scrX_1)$. 

A Shimura pair $(G,\scrX)$ of Hodge type is called of {\it PEL type}, if there exists an injective map $f:(G,\scrX)\hookrightarrow (\pmb{\text{GSp}}(W,\psi),\scrS)$ of Shimura pairs such that $G$ is the identity component of the intersection of $\pmb{\text{GSp}}(W,\psi)$ with the double centralizer of $G$ in $\pmb{GL}_W$. Here PEL stands for polarizations, endomorphisms, and level structures. In such a case, we say that  $f:(G,\scrX)\to (\pmb{\text{GSp}}(W,\psi),\scrS)$ is a {\it PEL type embedding}. If moreover $G$ itself is the intersection of $\pmb{\text{GSp}}(W,\psi)$ with the double centralizer of $G$ in $\pmb{\text{GL}}_W$, then we say that the Shimura pair $(G,\scrX)$ is of {\it moduli PEL type} and we say that $f:(G,\scrX)\to (\pmb{\text{GSp}}(W,\psi),\scrS)$ is a {\it moduli PEL type embedding}.

A Shimura pair $(G_1,\scrX_1)$ is called of {\it preabelian type} if there exists a Shimura pair $(G,\scrX)$ of Hodge type such that we have an isomorphism $(G^{\ad},\scrX^{\ad})\arrowsim (G_1^{\ad},\scrX_1^{\ad})$ of adjoint Shimura pairs. If moreover we can choose $(G,\scrX)$ such that the last isomorphism is defined by an isogeny $G^{\der}\to G_1^{\der}$ between the derived groups of $G$ and $G_1$, then we say that $(G_1,\scrX_1)$ is of {\it abelian type}.

We say $(G_1,\scrX_1)$ is {\it unitary}, if all simple factors of $G_{\dbC}^{\ad}$ are $\pmb{\text{PGL}}$ groups. Following [Va6, Def. 1.1] we say $(G_1,\scrX_1)$ {\it has compact factors}, if for each simple factor $G_0$ of $G_1^{\ad}$ there exists a simple factor of $G_{0,\dbR}$ which is compact. 

Let $\Sh(G_1,\scrX_1)$ be the canonical model of $(G_1,\scrX_1)$ over the reflex field $E(G_1,\scrX_1)$ of $(G_1,\scrX_1)$ (see [De1,2] and[Mi1,4]).  
The natural closed embedding of complex spaces $\scrX\times G(\dbA_f)\hookrightarrow \scrS\times \pmb{\text{GSp}}(W,\psi)(\dbA_f)$ induces naturally via passage to quotients a closed embedding $\Sh(G,\scrX)_{\dbC}\hookrightarrow \Sh(\pmb{\text{GSp}}(W,\psi),\scrS)_{\dbC}$. The reflex field $E(G,\scrX)$ is a number field which is the smallest subfield of $\dbC$ with the property that the last closed embedding is the pull-back of a closed embedding
$$\Sh(G,\scrX)\hookrightarrow \Sh(\pmb{\text{GSp}}(W,\psi),\scrS)_{E(G,\scrX)}\leqno (7)$$
(see [De1, Cor. 5.4]). Let $H:=K_p\cap G(\dbQ_p)$. As we have $\Sh(G,\scrX)(\dbC)=G(\dbQ)\backslash (\scrX\times G(\dbA_f))$ and $\Sh(\pmb{\text{GSp}}(W,\psi),\scrS)(\dbC)=\pmb{\text{GSp}}(W,\psi)(\dbQ)\backslash (\scrS\times \pmb{\text{GSp}}(W,\psi)(\dbA_f))$ (see [De2, Cor. 2.1.11]), it is easy to see that (7) induces naturally a closed embedding homomorphism
$$\Sh(G,\scrX)/H\hookrightarrow  \Sh(\pmb{\text{GSp}}(W,\psi),\scrS)_{E(G,\scrX)}/K_p.\leqno (8)$$
\noindent
{\bf 1.3. Integral canonical models.} Let $(G_1,\scrX_1)$ be a Shimura pair of abelian type such that the group $G_{1,\dbQ_p}$ is {\it unramified} (i.e., it has a Borel subgroup and it splits over  a finite, unramified extension of $\dbQ_p$). We recall that this equivalent to the fact $G_{1,\dbQ_p}$ extends to a reductive group scheme $G_{1,\dbZ_p}$ over $\dbZ_p$, cf. [Ti2]. Each compact, open subgroup of $G_1(\dbQ_p)=G_{1,\dbQ_p}(\dbQ_p)$ of the form $H_1:=G_{1,\dbZ_p}(\dbZ_p)$ is called {\it hyperspecial}. We refer to the triple $(G_1,\scrX_1,H_1)$ as a {\it Shimura triple of abelian type} (with respect to $p$). By a map $f:(G_1,\scrX_1,H_1)\to (\tilde G_1,\tilde\scrX_1,\tilde H_1)$ of Shimura triples of abelian type (with respect to $p$) we mean a map $f:(G_1,\scrX_1)\to (\tilde G_1,\tilde\scrX_1)$ of Shimura pairs such that the homomorphism $f(\dbQ_p):G_1(\dbQ_p)\to \tilde G_1(\dbQ_p)$ maps $H_1$ to $\tilde H_1$. 

As  the group $G_{1,\dbQ_p}$ is unramified, the reflex field $E(G_1,\scrX_1)$ is unramified above $p$ (cf. [Mi3, Cor. 4.7 (a)]). Thus the normalization $E(G_1,\scrX_1)_{(p)}$ of $\dbZ_{(p)}$ in $E(G_1,\scrX_1)$ is a finite, \'etale $\dbZ_{(p)}$-algebra. Let $\dbA_f^{(p)}$ be the ring of finite ad\`eles of $\dbQ$ with the $p$-component omitted; we have $\dbA_f=\dbA_f^{(p)}\times\dbQ_p$. 

\medskip\noindent
{\bf 1.3.1. Definitions.} {\bf (a)}  By an {\it integral model} of $\Sh(G_1,\scrX_1)/H_1$ over $E(G_1,\scrX_1)_{(p)}$ we mean a faithfully flat scheme $\scrN_1$ over $E(G_1,\scrX_1)_{(p)}$
together with a continuous right action of $G_1(\dbA_f^{(p)})$ on it in the sense of [De2, Subsubsect. 2.7.1], such that there exists a
$G_1(\dbA_f^{(p)})$-equivariant isomorphism
$$
\scrN_{1,E(G_1,\scrX_1)}\arrowsim\Sh(G_1,\scrX_1).
$$
The integral model $\scrN_1$ is said to be {\it smooth} (resp. {\it normal}) if there exists a compact, open
subgroup $\tilde H$ of $G_1(\dbA_f^{(p)})$ such that for every
inclusion $\tilde H_2\subseteq \tilde H_1$ of compact, open subgroups
of $\tilde H$, the natural morphism $\scrN_1/\tilde H_2\to \scrN_1/\tilde H_1$ is a finite, \'etale morphism
between smooth schemes (resp. between normal schemes) of finite type over $E(G_1,\scrX_1)_{(p)}$. 

We say that $\scrN_1$ is {\it quasi-projective} or {\it projective} if we can choose $\tilde H$ such that $\scrN_1/\tilde H$ is a quasi-projective or projective (respectively) $E(G_1,\scrX_1)_{(p)}$-scheme.

\smallskip
{\bf (b)} A regular, faithfully flat $E(G_1,\scrX_1)_{(p)}$-scheme $Y$ is called {\it healthy regular}, if for each open subscheme $U$ of $Y$ which contains $Y_{\dbQ}$ and all points of $Y$ of codimension $1$, every abelian scheme over $U$ extends uniquely to an abelian scheme over $Y$.

\smallskip
{\bf (c)} A scheme $Z$ over $E(G_1,\scrX_1)_{(p)}$ is said to have the {\it extension property} if for each  healthy
regular scheme $Y$ over $E(G_1,\scrX_1)_{(p)}$, 
every morphism $Y_{E(G_1,\scrX_1)}\to Z_{E(G_1,\scrX_1)}$ of $E(G_1,\scrX_1)$-schemes extends uniquely to a morphism $Y\to Z$ of $E(G_1,\scrX_1)_{(p)}$-schemes.

\smallskip
{\bf (d)}  A smooth integral model of $\Sh(G_1,\scrX_1)$ over $E(G_1,\scrX_1)_{(p)}$ that has the extension property is called an {\it integral canonical model} of $(G_1,\scrX_1,H_1)$ (or of $\Sh(G_1,\scrX_1)/H_1$ over $E(G_1,\scrX_1)_{(p)}$).

\medskip\noindent
{\bf 1.3.2. Two key facts.} Each regular, formally smooth scheme over $E(G,\scrX)_{(p)}$ is healthy regular, cf. [VZ, Cor. 5]. From this and Yoneda Lemma we get that (cf. [VZ, Cor. 30]): 

\medskip
{\bf (a)} an integral canonical model of $(G_1,\scrX_1,H_1)$ is uniquely determined up to a canonical isomorphism;

\smallskip
{\bf (b)} if we have a map $f:(G_1,\scrX_1,H_1)\to (\tilde G_1,\tilde\scrX_1,\tilde H_1)$ of Shimura triples of abelian type and if the integral canonical models $\scrN_1$ and $\tilde\scrN_1$ of $(G_1,\scrX_1,H_1)$ and $(\tilde G_1,\tilde\scrX_1,\tilde H_1)$ (respectively) exist, then the natural morphism $\Sh(G_1,\scrX_1)/H_1\to \Sh(\tilde G_1,\tilde\scrX_1)/\tilde H_1$ of $E(\tilde G_1,\tilde\scrX_1)$-schemes extends uniquely to a {\it functorial morphism} $\scrN_1\to\tilde\scrN_1$ of $E(\tilde G_1,\tilde\scrX_1)_{(p)}$-schemes. 

\smallskip
{\bf (c)} if the map $f:(G_1,\scrX_1,H_1)\to (\tilde G_1,\tilde\scrX_1,\tilde H_1)$ of (b) is injective, then similarly to (7) and (8), we get closed embeddings $\Sh(G_1,\scrX_1)\hookrightarrow \Sh(\tilde G_1,\tilde\scrX_1)_{E(\tilde G_1,\tilde\scrX_1)}$ and $\Sh(G_1,\scrX_1)/H_1\hookrightarrow \Sh(\tilde G_1,\tilde\scrX_1)_{E(\tilde G_1,\tilde\scrX_1)}/\tilde H_1$ (cf. also [Va1, Rm. 3.2.14] for the second one) and therefore the generic fibre of the functorial morphism 
$$\scrN_1\to\tilde\scrN_{1,E(G_1,\scrX_1)_{(p)}}\leqno (9)$$ 
is a closed embedding and one expects that (9) itself has nice properties as well.

\medskip
The below proposition is a $\dbZ_{(p)}$-variant of classical results of Satake and Deligne (see [Sa1,2] and [De2, Prop. 2.3.10]) and it is only slightly more general than [Va1, Thm. 6.5.1.1 a) to c)]; its proof is presented in Subsection 2.2. 
 
\bigskip\noindent
{\bf 1.4. Proposition.} {\it Let $(G_1,\scrX_1)$ be a simple, adjoint Shimura pair of abelian type. We assume that the group $G_{1,\dbQ_p}$ is unramified. Let $H_1$ be a hyperspecial subgroup of $G_1(\dbQ_p)$  i.e., be the group of $\dbZ_p$-valued points of a reductive group scheme $G_{1,\dbZ_p}$ over $\dbZ_p$ that extends $G_{1,\dbQ_p}$. Then there exists an injective map  $f:(G,\scrX)\hookrightarrow (\pmb{\text{GSp}}(W,\psi),\scrS)$ of Shimura pairs such that the following four properties hold:

\medskip
{\bf (i)} we have a natural identification $(G^{\ad},\scrX^{\ad})=(G_1,\scrX_1)$;

\smallskip
{\bf (ii)} there exists a $\dbZ$-lattice $L$ of $W$ which is self dual with respect to $\psi$ and for which the schematic closure $G_{\dbZ_{(p)}}$ of $G$ in $\pmb{\text{GSp}}(L\otimes_{\dbZ} \dbZ_{(p)},\psi)$ is a reductive group scheme;

\smallskip
{\bf (iii)} the natural quotient homomorphism $G_{\dbQ_p}\twoheadrightarrow G_{\dbQ_p}^{\ad}=G_{1,\dbQ_p}$ extends uniquely to a quotient homomorphism $G_{\dbZ_p}\twoheadrightarrow G_{\dbZ_p}^{\ad}=G_{1,\dbZ_p}$, where $G_{\dbZ_p}:=G_{\dbZ_{(p)}}\times_{\Spec\dbZ_{(p)}} \Spec\dbZ_p$;

\smallskip
{\bf (iv)} the semisimple group cover $G^{\der}$ of $G_1$ is the maximal one allowed by the abelian type (i.e., if $(G_2,\scrX_2)$ is any other Shimura pair of abelian type whose adjoint Shimura pair is $(G_1,\scrX_1)$, then the isogeny $G^{\der}\to G_1$ factors through the isogeny $G_2^{\der}\to G_1$).}

\medskip
The next three theorems are proved in Sections 3 to 5 (respectively). 

\bigskip\noindent
{\bf 1.5. Basic Theorem.} {\it Under the hypotheses of Proposition 1.4, we can choose the injective map of Shimura pairs $f:(G,\scrX)\hookrightarrow (\pmb{\text{GSp}}(W,\psi),\scrS)$ such that the properties 1.4 (i) to (iii) hold and moreover the following three additional properties hold as well:

\medskip

{\bf (i)} if $\scrN$ is the normalization of the schematic closure of $\Sh(G,\scrX)/H$ in $\scrM_{E(G,\scrX)_{(p)}}$ (this makes sense due to (8)), then $\scrN$ is the integral canonical model of $(G,\scrX,H)$ and is quasi-projective (here $H=G_{\dbZ_{(p)}}(\dbZ_p)$ is as in the end of Subsection 1.2); 

\smallskip
{\bf (ii)} the integral canonical model $\scrN_1$ of $(G_1,\scrX_1,H_1)$ exists and is quasi-projective;

\smallskip
{\bf (iii)} the functorial morphism $\scrN\to\scrN_1$ of $E(G_1,\scrX_1)_{(p)}$-schemes, is a pro-\'etale cover of an open closed subscheme of $\scrN_1$.}

\bigskip\noindent
{\bf 1.6. Main Theorem A.} {\it Let $(G_1,\scrX_1,H_1)$ be a Shimura triple of abelian type with respect to $p$. Then the following four properties hold:

\medskip
{\bf (a)} The integral canonical model $\scrN_1$ of $(G_1,\scrX_1,H_1)$ exists and it is quasi-projective.

\smallskip
{\bf (b)} Let $(G_1,\scrX_1,H_1)\to (G_2,\scrX_2,H_2)$ be a map of Shimura triples with respect to $p$ such that at the level of derived groups it induces an isogeny $G_1^{\der}\to G_2^{\der}$. The functorial morphism $\Sh(G_1,\scrX_1)/H_1\to\Sh(G_2,\scrX_2)/H_2$ of $E(G_2,\scrX_2)$-schemes extends uniquely to a morphism $f_1:\scrN_1\to\scrN_2$ of $E(G_2,\scrX_2)_{(p)}$-schemes, where $\scrN_2$ is the integral canonical model of $(G_2,\scrX_2,H_2)$. Then $f_1$ is a pro-\'etale cover of an open  closed subscheme of $\scrN_2$. 

\smallskip
{\bf (c)} Let $(G_1,\scrX_1,H_1)\hookrightarrow (G_2,\scrX_2,H_2)$ be an injective map of Shimura triples of abelian type. Let $\scrN_1\to\scrN_2$ be the functorial morphism of $E(G_2,\scrX_2)_{(p)}$-schemes, where $\scrN_2$ is the integral canonical model of $(G_2,\scrX_2,H_2)$. Then $\scrN_1$ is the normalization of the schematic closure of the closed subscheme $\Sh(G_1,\scrX_1)/H_1$ of $(\Sh(G_2,\scrX_2)/H_2)_{E(G_1,X_1)}$ in $\scrN_{2,E(G_1,X_1)_{(p)}}$. 

\smallskip
{\bf (d)} If the $\dbQ$--rank of the adjoint group $G_1^{\ad}$ is $0$, then $\scrN_1$ is projective.}

\bigskip\noindent
{\bf 1.7. Main Theorem B.} {\it We consider an injective map  $f:(G,\scrX)\hookrightarrow (\pmb{\text{GSp}}(W,\psi),\scrS)$ of Shimura pairs. Let $p$ be a prime such that there exists a $\dbZ$-lattice $L$ of $W$ with the properties that we have a perfect alternating form $\psi:L\times L\to\dbZ$ and the schematic closure $G_{\dbZ_{(p)}}$ of $G$ in $\pmb{\text{GL}}_{L\otimes_{\dbZ} \dbZ_{(p)}}$ is a reductive group scheme over $\dbZ_{(p)}$. Let $H=G_{\dbZ_{(p)}}(\dbZ_p)$ be as in the end of Subsection 1.2. Then the following two properties hold:

\medskip
{\bf (a)} If $\scrN$ is the normalization of the schematic closure of $\Sh(G,\scrX)/H$ in $\scrM_{E(G,\scrX)_{(p)}}$ (this makes sense due to (8)), then $\scrN$ is the integral canonical model of $(G,\scrX,H)$ and is quasi-projective. 

\smallskip
{\bf (b)} If $O_{(v)}$ is the localization of $E(G,\scrX)_{(p)}$ at a maximal ideal $v$ and if $\scrN^{\m}_{O_{(v)}}$ is the $G(\dbA_f^{(p)})$-invariant, open subscheme of $\scrN_{O_{(v)}}$ defined in [Va15, Subsubsect. 3.5.1] (the definition is recalled in Section 5), then we have $\scrN^{\m}_{O_{(v)}}=\scrN_{O_{(v)}}$ (i.e., the motivic conjecture of Milne mentioned in [Va11] and [Mi5, Sect. 5] holds).}

\bigskip\noindent
{\bf 1.8. Previous results pertaining to Theorems 1.5 to 1.7.} They can be grouped as follows.

\medskip
{\bf (i)}  Mumford proved the existence of integral canonical models of Siegel modular varieties. More precisely, the $\dbZ_{(p)}$-scheme $\scrM$ together with the natural action of $\pmb{\text{GSp}}(W,\psi)(\dbA_f^{(p)})$ on it, is an integral canonical model of $(\pmb{\text{GSp}}(W,\psi),\scrS,K_p)$. Artin's method can be used to regain this result (see [Ar1,2], and [FC, Ch. I, Subsect. 4.11]).

\smallskip
{\bf (ii)} If $f:(G,\scrX)\hookrightarrow (\pmb{\text{GSp}}(W,\psi),\scrS)$ is a moduli PEL type embedding and if the condition 1.4 (ii) holds, then Zink proved that the schematic closure of $\Sh(G,\scrX)/H$ in $\scrM_{E(G,\scrX)_{(p)}}$ (this makes sense due to (8)) is the integral canonical model of $(G,\scrX,H)$ (see [Zi]). This result was reobtained in [LR]. 

\smallskip
{\bf (iii)} If $f:(G,\scrX)\hookrightarrow (\pmb{\text{GSp}}(W,\psi),\scrS)$ is a PEL type embedding, if the condition 1.4 (ii) holds, and if $p>2$, then Kottwitz pointed out that the arguments of [LR] can be used to get as well that the schematic closure of $\Sh(G,\scrX)/H$ in $\scrM_{E(G,\scrX)_{(p)}}$ is the integral canonical model of $(G,\scrX,H)$ (see [Ko]).

\smallskip
{\bf (iv)} In [Va1] it is proved that Theorem 1.5, Theorem 1.6 (a) to (c), Theorem 1.7 (a), and a weaker form of Theorem 1.7 (b) hold provided $p\ge 5$.  

\smallskip
{\bf (v)} See [VZ, Cors. 5 and 30] and [Va7, Thm. 1.3] for two corrections to [Va1] in connection to (iv). More precisely:

\medskip\noindent
$\bullet$ the original argument of Faltings in [Va1, Subsubsect. 3.2.17, Step B, last paragraph] and of Faltings and Chai in [FC, top of p. 184] were incorrect and they have been corrected by [Va3, Prop. 4.1] and by [VZ, Sect. 5] (respectively);

\smallskip\noindent
$\bullet$ the proof of Theorem 1.6 (a) for $p\ge 5$ and for the cases when $G_{1,\dbC}^{\ad}$ has simple factors isomorphic to $\pmb{\text{PGL}}_{pm}$ for some $m\in\dbN^{\ast}$ was partially incorrect in [Va1]; this has been corrected by [Va7, Thm. 1.3] (cf. [Va7, Appendix, E.3]).

\smallskip
{\bf (vi)} In [Va7] it is proved that Theorems 1.5 and 1.6 (a) hold provided $(G_1,\scrX_1)$ is a unitary Shimura pair. 

\smallskip
{\bf (vii)} Theorem 1.6 (d) is only a direct consequence of the previous results [Mo], [Pa], [Va6,7], and [Lee].

\smallskip
{\bf (viii)} In [Va12,13] it is shown that Kottwitz's result (see (iii)) holds even if $p=2$. 

\smallskip
{\bf (ix)} In [Va15] it is proved that Theorems 1.5 and 1.7 hold if the adjoint Shimura pair $(G^{\ad},\scrX^{\ad})$ has compact factors.

\smallskip
{\bf (x)} In [Ki3] it is claimed that Theorems 1.5, 1.6 (a), and 1.7 hold if either $p>2$ or $p=2$ and a very technical condition holds (strictly speaking, [Ki3] works with a weaker form of Proposition 1.4 (ii)). The paper [Ki3] does not bring any new conceptual ideas to [Va1,7,11,12] (being in fact only a variation of the ideas of loc. cit.). This variation was made possible due to advances in the theory of crystalline representations achieved by Fontaine, Breuil, Berger, and Kisin (see [Ki1,2], etc.).

\bigskip\noindent
{\bf 1.9. On the strategy to prove 1.5 and 1.6.} Main Theorems A and B are direct consequences of (the proof of) Theorem 1.5 and of the methods developed in [Va1--13]. Thus we will detail here only on the new strategy to prove Theorem 1.5. 

It is well known that $(G_1,\scrX_1)$ is of one of the following five disjoint types: $A_n$ (with $n\ge 1$), $B_n$ (with $n\ge 3$), $C_n$ (with $n\ge 2$), $D_n^{\dbH}$ (with $n\ge 4$), and $D_n^{\dbR}$ (with $n\ge 4$). These types were introduced in [De2, Table 2.3.8]. For instance, $(G_1,\scrX_1)$ is of $A_n$, $B_n$, or $C_n$ type if and only if all simple factors of $G_{1,\dbC}$ are of $A_n$, $B_n$, or $C_n$ (respectively) Lie type. To prove Theorem 1.5, we considers three disjoint cases:

\medskip
{\bf (PELNONCOMP)} all simple factors of $G_{1,\dbR}$ are non-compact and $(G_1,\scrX_1)$ is of either $A_n$ (with $n\Ge 1$) or $C_n$ (with $n\Ge 2$) or $D_n^{\dbH}$ (with $n\Ge 4$) type;

\smallskip
{\bf (COMP)} there exists a simple factor of $G_{1,\dbR}$ which is compact (i.e., the Shimura pair $(G_1,\scrX_1)$ has compact factors);

\smallskip
{\bf (SPINNONCOMP)} all simple factors of $G_{1,\dbR}$ are non-compact and $(G_1,\scrX_1)$ is of either $B_n$ (with $n\Ge 3$) or $D_n^{\dbR}$ (with $n\Ge 4$) type.

\medskip
The proof of Theorem 1.5 in the (PELNONCOMP) case relies on the results 1.8 (ii), (iii), and (viii) and its essence is very well documented in the literature (see [Va7] for the $A_n$ type; the $C_n$ and $D_n^{\dbH}$ type cases are similar).

The proof of Theorems 1.5 (i) in the (COMP) case was presented in [Va15]. 

The proof of Theorem 1.5 (i) in the (SPINNONCOMP) case involves four basic results (see Subsection 3.6); we list them here using bullets. 

\medskip
$\bullet$ It is well known that $E(G_1,\scrX_1)=\dbQ$ and that moreover one can choose the injective map  $f:(G,\scrX)\hookrightarrow (\pmb{\text{GSp}}(W,\psi),\scrS)$ such that we also have $E(G,\scrX)=\dbQ$ (see [De2]). 

\smallskip
$\bullet$ As $E(G,\scrX)=\dbQ$, a standard versal argument involving $F$-isocrystals shows that we can choose the injective map  $f:(G,\scrX)\hookrightarrow (\pmb{\text{GSp}}(W,\psi),\scrS)$ such that the {\it ordinary locus} of the special fibre $\scrL$ of $\scrN$ is Zariski dense in $\scrL$ (cf. Proposition 3.6.1). 

\smallskip
$\bullet$ It is well known that the ordinary points of $\scrN$ belong to the regular, formally smooth locus of $\scrN$ over $\dbZ_{(p)}$ (see [No]). 

\smallskip
$\bullet$ Based on [Va15, Thm. 1.7 (a) and (b)], the intersection of the smooth locus of $\scrN$ with $\scrL$ is an open closed subscheme of $\scrL$ and thus (due to the last two bullets) it is $\scrL$ itself.

\medskip
The passage from Theorem 1.5 (i) to Theorem 1.5 (ii) and (iii) in the (COMP) and (SPINNONCOMP)  cases is the same as in [Va1] if $p>2$ and it is very much the same as in [Va7] if $p=2$. More precisely, regardless of what the prime $p$ is we show that we can choose the injective map  $f:(G,\scrX)\hookrightarrow (\pmb{\text{GSp}}(W,\psi),\scrS)$ such that it factors through an injective map $f_2:(G_2,\scrX_2)\hookrightarrow (\pmb{\text{GSp}}(W,\psi),\scrS)$ with the properties that: (a) $(G_2^{\ad},\scrX_2^{\ad})$ is an adjoint, unitary Shimura pair, and (b) the resulting injective map $(G,\scrX)\hookrightarrow (G_2,\scrX_2)$ induces an injective map $(G_1,\scrX_1)\hookrightarrow (G_2^{\ad},\scrX_2^{\ad})$ at the level of adjoint Shimura pairs. Using such a factorization we conclude that Theorem 1.5 (ii) and (iii) follow from Theorem 1.5 (i) and [Va7].

\bigskip
\noindent
{\boldsectionfont 2. Preliminaries}
\bigskip
 
Subsection 2.1 recalls standard notations on reductive group schemes and modules. Subsection 2.2 proves Proposition 1.4. Subsection 2.3 recalls the notion of a cover between Shimura triples of abelian type introduced in [Va7, Subsect. 2.4]. Proposition 2.4 is the very essence of Theorem 1.6 (c).

\bigskip\noindent
{\bf 2.1. Notations.} A group scheme $\scrH$ over an affine scheme $\Spec R$ is called {\it reductive}, if it is smooth and affine and its fibres are connected and have trivial unipotent radicals (cf. [DG, Vol. III, Exp. XIX, Def. 2.7]). Let $\scrH^{\der}$, $Z(\scrH)$, $\scrH^{\ab}$, and $\scrH^{\ad}$ be the {\it derived group scheme} of $\scrH$, the {\it center} of $\scrH$, the {\it maximal commutative quotient} of $\scrH$, and the {\it adjoint group scheme} of $\scrH$ (respectively). Thus we have $\scrH^{\ad}=\scrH/Z(\scrH)$ and $\scrH^{\ab}=\scrH/\scrH^{\der}$, cf. [DG, Vol. III, Exp. XXII, Def. 4.3.6 and Thm. 6.2.1]. For a finite, \'etale morphism $\Spec R\to \Spec R_0$, let $\Res_{R/R_0} \scrH$ be the reductive group scheme over $\Spec R_0$ obtained from $\scrH$ through the {\it Weil restriction of scalars} (see [BLR, Ch. 7, Sect. 7.6] and [Va2, Subsect. 2.3]). 

If $M$ is a free module of finite rank over a  commutative $\dbZ$-algebra $R$, then let $M^{\vee}:=\Hom_R(M,R)$ and let $\scrT(M):=\oplus_{s,t\in\dbN} M^{\otimes s}\otimes_R M^{{\vee}\otimes t}$.

\bigskip\noindent
{\bf 2.2. Proof of 1.4.} Let $(G_1,\scrX_1)$ be a simple, adjoint Shimura pair of abelian type. We assume that the group $G_{1,\dbQ_p}$ is unramified. Let $H_1$ be a hyperspecial subgroup of $G_1(\dbQ_p)$  i.e., the group of $\dbZ_p$-valued points of a reductive group scheme $G_{1,\dbZ_p}$ over $\dbZ_p$ that extends $G_{1,\dbQ_p}$. In this subsection we prove Proposition 1.4. Let $F_1$ be a number field such that we have an isomorphism $G_1\arrowsim \text{Res}_{F_1/\dbQ} J_1$, where $J_1$ is an absolutely simple adjoint group over $F_1$ (see [Ti1, Subsubsect. 3.1.2]). The number field $F_1$ is uniquely determined up to $\Gal(\Bbb Q)$-conjugation (i.e., up to isomorphism).

As $G_{1,\dbR}$ is an inner form of its compact form (cf. [De2, p. 255]), it is a product of absolutely simple, adjoint groups over $\dbR$. This implies that the number field $F_1$ is totally real. As $G_{1,\dbQ_p}$ splits over an unramified extension $F_p$ of $\dbQ_p$, the $F_p$-algebra $F_1\otimes_{\dbQ} F_p$ is isomorphic to $F_p^{[F:\dbQ]}$. Therefore $F_1$ is unramified over $p$. 

Let $E_1$ be a totally imaginary quadratic extension of $F_1$ unramified over $p$. If $(G_1,\scrX_1)$ is of $A_n$ type with $n\ge 2$, then we take $E_1$ to be uniquely determined by the property that $\Gal(E_1)$ acts trivially on the  ending notes of the Dynkin diagram of $G_{1,\overline{\dbQ}}$ (see either [De2] or [Va8, Subsect. 2.2] for the action of $\Gal(\dbQ)$ on the Dynkin diagram of $G_{1,\overline{\dbQ}}$; note that $E_1$ is unramified over $p$ as $G_1$ has a maximal torus which splits over a finite Galois extension of $\dbQ$ unramified over $p$). We consider an injective map $f:(G,\scrX)\hookrightarrow (\pmb{\text{GSp}}(W,\psi),\scrS)$ of Shimura pairs such that the properties 1.4 (i) and (iv) hold and moreover the torus $Z^0(G)$ is a subtorus of $\text{Res}_{E_1/\dbQ} \dbG_m$, cf. [De2, Prop. 2.3.10] and its proof if $(G_1,\scrX_1)$ is not of $A_n$ type with $n\ge 2$ and cf. the modifications to loc. cit. made in [Va7, Prop. 3.2] if $(G_1,\scrX_1)$ is of  $A_n$ type with $n\ge 2$. If $(G_1,\scrX_1)$ is of either $D_n^{\dbR}$ or $D_n^{\dbH}$ type with $n\ge 4$, then the fact that the torus $Z^0(G)$ is a subtorus of $\text{Res}_{E_1/\dbQ} \dbG_m$ follows as well from either [De2, Rm. 2.3.13] or (proof of) [Va8, Thm. 4.8]. We conclude that the torus $Z^0(G)$ splits over a Galois extension of $\dbQ$ unramified above $p$. Thus the torus $Z^0(G)_{\dbQ_p}$ is unramified. As $G_{\dbQ_p}$ is isogeneous to  $Z^0(G)_{\dbQ_p}\times_{\dbQ_p} G_{1,\dbQ_p}$, it is also unramified.

Let $H$ be a hyperspecial subgroup of $G(\dbQ_p)$. Let $G_{\dbZ_{(p)}}$ be the unique reductive group scheme over $\dbZ_{(p)}$ whose generic fibre is $G$ and whose group of $\dbZ_p$-valued points is $H$, cf. [Va7, Lem. 2.3 (a)]. Let $\tilde H_1:=G^{\ad}_{\dbZ_p}(\dbZ_p)$; it is a hyperspecial subgroup of $G(\dbQ_p)$. The hyperspecial subgroups of $G_1(\dbQ_p)$ are $G_1(\dbQ_p)$-conjugate, cf. [Ti2, p. 47]. Thus there exists an element $g\in G_1(\dbQ_p)$ such that we have $\tilde H_1=gH_1 g^{-1}$. By replacing $H$ with $g^{-1}H g$, $\tilde H_1$ gets replaced by $H_1=g^{-1}\tilde H_1 g$. Thus we can assume that $\tilde H_1=H_1$. This implies that $G_{1,\dbZ_p}=G_{\dbZ_p}^{\ad}$, cf. loc. cit. or [Va7, Lem. 2.3 (a)]. Thus the property 1.4 (iii) holds.

Let $L$ be a $\dbZ$-lattice of $W$ which is self-dual with respect to $\psi$ (i.e., $\psi$ induces a perfect alternating form $\psi:L\times L\to \dbZ$). From [Va15, Lem. 4.2.1] we get that we can modify the $\dbZ$-lattice $L$ and the injective map $f:(G,\scrX)\hookrightarrow (\pmb{\text{GSp}}(W,\psi),\scrS)$ of Shimura pairs such that moreover $L_{(p)}:=L\otimes_{\dbZ} \dbZ_{(p)}$ is a $G_{\dbZ_{(p)}}$-module. The resulting homomorphism $G_{\dbZ_{(p)}}\to \pmb{\text{GL}}_{L_{(p)}}$ is a closed embedding, cf. [Va4, Thm. 1.1] and [Va15, Fact 2.3.1]. Thus the property 1.4 (ii) holds as well.\endproof

\medskip\noindent
{\bf 2.2.1. Notations related to 1.4 and 1.7.} We consider an injective map  $f:(G,\scrX)\hookrightarrow (\pmb{\text{GSp}}(W,\psi),\scrS)$ of Shimura pairs. Let $p$ be a prime such that there exists a $\dbZ$-lattice $L$ with the properties that we have a perfect alternating form $\psi:L\times L\to\dbZ$ and the schematic closure $G_{\dbZ_{(p)}}$ of $G$ in $\pmb{\text{GL}}_{L_{(p)}}$ is a reductive group scheme over $\dbZ_{(p)}$; here $L_{(p)}:=L\otimes_{\dbZ} \dbZ_{(p)}$. Let $\psi^{\vee}$ be the perfect alternating form on $L_{(p)}^{\vee}$ that is defined naturally by $\psi$. Let $(v_{\alpha})_{\alpha\in\scrJ}$ be a family of tensors of $\scrT(W^{\vee})$ such that $G$ is the subgroup of $\pmb{\text{GL}}_{W^{\vee}}$ that fixes $v_{\alpha}$ for all $\alpha\in\scrJ$, cf. [De3, Prop. 3.1.2 c)]. Let $H=G_{\dbZ_{(p)}}(\dbZ_p)=G(\dbQ_p)\cap K_p$ be as in the end of Subsection 1.2.

Let $\scrN$ be the normalization of the schematic closure of $\Sh(G,\scrX)/H$ in $\scrM_{E(G,\scrX)_{(p)}}$. Let $\scrN^{\s}$ be the open subscheme of $\scrN$ which is the formally smooth locus of $\scrN$ over $E(G,\scrX)_{(p)}$. We have an identity $\scrN^{\s}_{E(G,\scrX)}=\scrN_{E(G,\scrX)}$, cf. [Va15, Lem. 2.2.2]. Then $\scrN$ is a quasi-projective integral canonical model of $(G,\scrX,H)$ (i.e., property 1.5 (i) holds) if and only if we have $\scrN^{\s}=\scrN$, cf. [Va1, Cor. 3.4.4]. Let $(\scrA,\Lambda_{\scrA})$ be the pull-back to $\scrN$ of the universal principally polarized abelian scheme over $\scrM$.

Let $k$ be an algebraically closed field of characteristic $p$ and countable transcendental degree. Let $W(k)$ be the ring of Witt vectors with coefficients in $k$ and let $B(k)=W(k)[{1\over p}]$ be its field of fractions. Let $\sigma:=\sigma_k$ be the Frobenius automorphism of $k$, $W(k)$, and $B(k)$.

\bigskip\noindent
{\bf 2.3. Definition.} A morphism  $f:(G_1,\scrX_1,H_1)\to (\tilde G_1,\tilde\scrX_1,\tilde H_1)$ between Shimura triples of abelian type is called a {\it cover}, if the following two properties hold:

\medskip
{\bf (i)} the group $G_1$ surjects onto $\tilde G_1$, and

\smallskip
{\bf (ii)} the kernel $\Ker(f)$ is a torus of $Z(G_1)$ with the property that for every field $K$ of characteristic $0$, the cohomology group $H^1(K,\Ker(f)_K)$ is trivial. 

\medskip
Each cover $f:(G_1,\scrX_1,H_1)\to (\tilde G_1,\tilde\scrX_1,\tilde H_1)$ induces at the level of adjoint triples an isomorphism $f^{\ad}:(G^{\ad}_1,\scrX_1^{\ad},H_1^{\ad})\arrowsim (\tilde G_1^{\ad},\tilde\scrX_1^{\ad},\tilde H_1^{\ad})$.

\bigskip\noindent
{\bf 2.4. Proposition.} {\it Let $(G_1,\scrX_1,H_1)\hookrightarrow (G_2,\scrX_2,H_2)$ be an injective map of Shimura triples of abelian type with respect to $p$. For $j\in\{1,2\}$ we assume that the integral canonical model of $(G_j,\scrX_j,H_j)$ exists and is quasi-projective. We view $\scrN_1$ as a $\scrN_2$-scheme via the functorial morphism $\scrN_1\to\scrN_2$ of $E(G_2,\scrX_2)_{(p)}$-schemes. Then $\scrN_1$ is the normalization $\scrP_1$ of the schematic closure of $\Sh(G_1,\scrX_1)/H_1$ in $\scrN_{2,E(G_1,\scrX_1)_{(p)}}$ (we recall from the fact 1.3.2 (c) that the functorial morphism $\Sh(G_1,\scrX_1)/H_1\to \Sh(G_2,\scrX_2)_{E(G_1,\scrX_1)}/H_2$ is a closed embedding).} 

\medskip
\proof
It is known that $\scrP_1$ is a normal integral model of $(G_1,\scrX_1,H_1)$ that has the extension property, cf. [Va1, Prop. 3.4.1]. We have a natural morphism $a:\scrN_1\to\scrP_1$ of $E(G_1,\scrX_1)_{(p)}$-schemes whose generic fibre is the identity automorphism of $\Sh(G_1,\scrX_1)/H_1$. The morphism $a$ is a pro-\'etale cover of a morphism $a_{H_0}:\scrN_1/H_0\to\scrP_1/H_0$ of normal $E(G_1,\scrX_1)_{(p)}$-schemes of finite type, where $H_0$ is a small enough compact, open subgroup of $G_1(\dbA_f^{(p)})$ (cf. Definition 1.3.1 (a)). As $\scrN_2$ is quasi-projective and as $E(G_1,\scrX_1)_{(p)}$ is an excellent ring, it is easy to see that the $E(G_1,\scrX_1)_{(p)}$-scheme $\scrP_1/H_0$ is quasi-projective. Thus $a_{H_0}$ is a quasi-projective morphism between normal, flat $E(G_1,\scrX_1)_{(p)}$-schemes of finite type whose generic fibre is an isomorphism. As each discrete valuation ring of mixed characteristic $(0,p)$ is a healthy regular scheme, the morphism $a$ satisfies the valuative criterion of properness with respect to such discrete valuation rings. From the last two sentences we get that $a_{H_0}$ is in fact a projective morphism. 

We consider the largest open subscheme $\scrO_1$ of $\scrP_1$ which contains $\Sh(G_1,\scrX_1)/H_1$ and for which the morphism $a^{-1}(\scrO_1)\to\scrO_1$ is an isomorphism. As $\scrN_1$ has the extension property and is a smooth integral model (cf. Definition 1.3.1 (d)) and as each regular, formally smooth scheme over $E(G_1,\scrX_1)_{(p)}$ is healthy, we get that $\scrO_1$ is in fact the formally smooth locus of $\scrP_1$ over $E(G_1,\scrX_1)_{(p)}$. Obviously $\scrO_1$ is $H_0$-invariant. Thus the projective morphism $a_{H_0}:\scrN_1/H_0\to\scrP_1/H_0$ is an isomorphism above $\scrO_1/H_0$. As $a_{H_0}$ is projective, we have an inequality $\text{codim}_{\scrP_1/H_0}((\scrP_1/H_0)\setminus (\scrO_1/H_0))\ge 2$. 

Let $\scrY$ be the set of points of $\scrN_1/H_0$ which are of codimension $1$ and which do not belong to $a_{H_0}^{-1}(\scrO_1/H_0)\arrowsim\scrO_1/H_0$. To prove the proposition it suffices to show that $a_{H_0}$ is an isomorphism. To check that $a_{H_0}$ is an isomorphism, it suffices to show that the set $\scrY$ is empty (this is so as the projective morphism $a_{H_0}$ is a blowing up of a closed subscheme of $\scrP_1/H_0$; the proof of this is similar to [Ha, Ch. II, Thm. 7.17]). 

We show that the assumption that the set $\scrY$ is non-empty leads to a contradiction. Let $\scrC$ be the open subscheme of $\scrN_1/H_0$ which contains $\scrN_{1,E(G_1,\scrX_1)}/H_0$ and for which $\scrE_1:=\scrC_{\dbF_p}$ is the union of all connected components of special fibres of $\scrN_1/H_0$ whose generic points belong to $\scrY$. Let $d_1:=\dim(\scrE_1)$. The image $\scrE_2:=a_{H_0}(\scrE_1)$ has dimension less than $d_1$ and is contained in the non-smooth locus of $\scrP_1/H_0$. The morphism $\scrC\to\scrP_1/H_0$ factors through the dilatation $\scrV$ of $\scrP_1/H_0$ centered on the reduced scheme of the non-smooth locus of $\scrP_1/H_0$, cf. the universal property of dilatations (see [BLR, Ch. 3, 3.2, Prop. 3.1 (b)]). But $\scrV$ is an affine $\scrP_1/H_0$-scheme and thus the image of the projective $\scrP_1/H_0$-scheme $\scrE_1$ in $\scrV$ has the same dimension as $\scrE_2$. By repeating the process we get that the image of $\scrE_1$ in a smoothening $\scrV_{\infty}$ of $\scrP_1/H_0$ obtaining via a sequence of dilatations centered on non-smooth loci (see [BLR, Ch. 3, Thm. 3 of Sect. 3.1 and Thm. 2 of Sect. 3.4]), has dimension $\dim(\scrE_2)$ and thus it has dimension less than $d_1$. But each discrete valuation ring of $\scrV_{\infty}$ dominates a local ring of $\scrN_1/H_0$ (as $a_{H_0}$ is a projective morphism) and therefore (due to the existence of the morphism $\scrC\to \scrV_{\infty}$) it is also a local ring of 
$\scrN_1/H_0$. As $\scrV_{\infty}$ has at least one discrete valuation ring of mixed characteristic $(0,p)$ which is not a local ring of $\scrO_1/H_0$, we get that this discrete valuation ring is the local ring of a point in $\scrY$. Thus $\text{Im}(\scrE_1\to\scrV_{\infty})$ has dimension $d_1$. Contradiction. Thus $\scrY=\emptyset$ and $a:\scrN_1\to\scrP_1$ is an isomorphism.\endproof

\bigskip
\noindent
{\boldsectionfont 3. The proof of Theorem 1.5}
\bigskip

Subsections 3.1 to 3.6 present the six cases that are required to prove Theorem 1.5 (i). Theorems 1.5 (ii) and (iii) are proved in Subsection 3.8 based on Proposition 3.7. Corollary 3.9 is the very essence of the proof of Theorem 1.6 (a) and (b).

Let $(G_1,\scrX_1)$ be a simple, adjoint Shimura pair of abelian type. We assume that the group $G_{1,\dbQ_p}$ is unramified. Let $H_1$ be a hyperspecial subgroup of $G_1(\dbQ_p)$  i.e., the group of $\dbZ_p$-valued points of a reductive group scheme $G_{1,\dbZ_p}$ over $\dbZ_p$ that extends $G_{1,\dbQ_p}$. Until Subsection 1.3 we will consider an injective map $f:(G,\scrX)\hookrightarrow (\pmb{\text{GSp}}(W,\psi),\scrS)$ of Shimura pairs such that the properties 1.4 (i) to (iv) hold. Let $L_{(p)}:=L\otimes_{\dbZ} \dbZ_{(p)}$ and let $\scrO:=\End_{\dbZ_{(p)}}(L_{(p)})\cap\{e\in \End_{\dbQ}(W)|e\;\text{is}\;\text{fixed}\;\text{by}\;G\}$.

\bigskip\noindent
{\bf 3.1. Case 1: the unitary case.}
We assume that the Shimura pair $(G_1,\scrX_1)$ is unitary. Based on [Va7, Prop. 3.2], we can assume that the $\dbZ_{(p)}$-algebra $\scrO$ is semisimple and that $G_{\dbZ_{(p)}}$ is the subgroup scheme of $\pmb{\text{GSp}}(L_{(p)},\psi)$ that fixes $\scrO$. From this and [Va7, Cor. 4.1.1] we get that the property 1.5 (i) holds. The fact that the properties 1.5 (ii) and (iii) hold follows from [Va7, Thm. 5.1 (a)]  and [Va7, Thm. 5.1 (b)] (respectively). Thus Theorem 1.5 holds if the Shimura pair $(G_1,\scrX_1)$ is unitary.

\bigskip\noindent
{\bf 3.2. Case 2: the totally non-compact $C_n$ type case.} We assume that the Shimura pair $(G_1,\scrX_1)$ is of $C_n$ type and that all simple factors of $G_{1,\dbR}$ are non-compact. Then arguments entirely similar to the ones of [Va7, Prop. 3.2] show that we can assume that the $\dbZ_{(p)}$-algebra $\scrO$ is semisimple and that $G_{\dbZ_{(p)}}$ is the subgroup scheme of $\pmb{\text{GSp}}(L_{(p)},\psi)$ that fixes $\scrO$. From this and [Zi, Subsect. 3.5], [LR], and [Ko, Sect. 5] we get that the property 1.5 (i) holds and in fact $\scrN$ is the schematic closure of $\Sh(G,\scrX)/H$ in $\scrM_{E(G,\scrX)_{(p)}}$. Using this, [Va7, Subsects. 4.1 to 4.3 and 5.1] can be adapted to show that the properties 1.5 (ii) and (iii) hold as well. Thus Theorem 1.5 holds if the Shimura pair $(G_1,\scrX_1)$ is of $C_n$ type and all simple factors of $G_{1,\dbR}$ are non-compact.

\bigskip\noindent
{\bf 3.3. Case 3: the totally non-compact $D_n^{\dbH}$ type case with $p>2$.}
We assume that $p>2$, that the Shimura pair $(G_1,\scrX_1)$ is of $D_n^{\dbH}$ type, and that all simple factors of $G_{1,\dbR}$ are non-compact. Then arguments similar to the ones of [Va7, Prop. 3.2] show that we can assume that the $\dbZ_{(p)}$-algebra $\scrO$ is semisimple and that $G_{\dbZ_{(p)}}$ is the identity component of the subgroup scheme of $\pmb{\text{GSp}}(L_{(p)},\psi)$ that fixes $\scrO$. From this and [Ko, Sect. 5] we get that the property 1.5 (i) holds and $\scrN$ is the schematic closure of $\Sh(G,\scrX)/H$ in $\scrM_{E(G,\scrX)_{(p)}}$ and is quasi-projective. Using this, [Va7, Subsects. 4.1 to 4.3 and 5.1] can be adapted to show that the properties 1.5 (ii) and (iii) hold as well. Thus Theorem 1.5 holds if $p>2$, the Shimura pair $(G_1,\scrX_1)$ is of $D_n^{\dbH}$ type, and all simple factors of $G_{1,\dbR}$ are non-compact.

For the sake of completeness we present a second way to argue that properties 1.5 (ii) and (iii) hold. Let $(G^\prime,\scrX^\prime,H^\prime)\to (G_1,\scrX_1,H_1)$ be a map of Shimura triples for which the following two properties hold (cf. [MS, 3.4] and [Va1, Rm. 3.2.7 10)]): (i) it is a cover in the sense of Subsection 2.3 and (ii) we have identities $E(G^\prime,\scrX^\prime)=E(G_1,\scrX_1)$ and $G^{\prime,\der}=G^{\der}$. As $G^{\prime,\der}=G^{\der}$, the integral canonical model $\scrN^\prime$ of $(G^\prime,\scrX^\prime,H^\prime)$ over $E(G^\prime,\scrX^\prime)_{(p)}=E(G_1,\scrX_1)_{(p)}$ exists (cf. [Va7, Prop. 4.2.3 (a)] and the fact that the integral canonical model $\scrN$ of $(G,\scrX,H)$ exists). As $\scrN$ is quasi-projective, from [Va7, Prop. 2.4.3 (c)] we get that $\scrN^\prime$ is also quasi-projective. The order of the center of the simply connected semisimple group cover of $G_1$ is a power of $2$ and thus it is relative prime to $p$. From the last two sentences and [Va1, Thm. 6.2 (a)] we get that the integral canonical model $\scrN_1$ of $(G_1,\scrX_1,H_1)$ exists and that the functorial morphism $\scrN^\prime\to\scrN_1$ is a pro-\'etale cover of an open closed subscheme of $\scrN_1$. From this and [Va7, Prop. 4.2.3 (c)] we get that the functorial morphism $\scrN^\prime\to\scrN_1$ is a pro-\'etale cover of an open closed subscheme of $\scrN_1$. Thus properties 1.5 (ii) and (iii) hold.

\bigskip\noindent
{\bf 3.4. Case 4: the totally non-compact $D_n^{\dbH}$ type case with $p=2$.}
We assume that $p=2$, that the Shimura pair $(G_1,\scrX_1)$ is of $D_n^{\dbH}$ type, and that all simple factors of $G_{1,\dbR}$ are non-compact. Then arguments similar to the ones of [Va7, Prop. 3.2] show that we can assume that the $\dbZ_{(p)}$-algebra $\scrO$ is semisimple and that $G_{\dbZ_{(p)}}$ is the schematic closure in $\pmb{\text{GSp}}(L_{(p)},\psi)$ of the identity component of the subgroup scheme of $\pmb{\text{GSp}}(W,\psi)$ that fixes $\scrO$. From this and [Va14, Thm. 1.3] we get that the property 1.5 (i) holds.

\bigskip\noindent
{\bf 3.5. Case 5: the non-unitary, compact factors case.}
We assume that the Shimura pair $(G_1,\scrX_1)$ has compact factors and is not unitary. From [Va15, Thm. 1.7 (b)] we get that the property 1.5 (i) holds. As  $(G_1,\scrX_1)$ is not unitary, it is of $B_n$, $C_n$, $D_n^{\dbH}$, or $D_n^{\dbR}$ type. Thus the order of the center of the simply connected semisimple group cover of $G_1$ is a power of $2$. Thus, if moreover $p>2$, then as in the second paragraph of Section 3.3 we argue that properties 1.5 (ii) and (iii) also hold.

\bigskip\noindent
{\bf 3.6. Case 6: the totally non-compact $B_n$ and $D_n^{\dbR}$ types case.}
We assume that the Shimura pair $(G_1,\scrX_1)$ is of either $B_n$ or $D_n^{\dbR}$ type and that all simple factors of $G_{1,\dbR}$ are non-compact. We have $E(G_1,\scrX_1)=\dbQ$, cf. [De2, Rm. 2.3.12].  We can assume that $Z^0(G)=\dbG_m$ and $E(G,\scrX)=\dbQ$, cf. [De2, Rm. 2.3.13]. This implies that there exists a cocharacter $\mu_0:\dbG_m\to G_{\dbZ_p}$ whose extension to $\dbC$ is $G(\dbC)$-conjugate to the Hodge cocharacters $\mu_x:\dbG_m\to G_{\dbC}$ associated naturally to $x\in\scrX$, cf. [Mi3, Prop. 4.6]. We have a direct sum decomposition $L_p^{\vee}:=L^{\vee}\otimes_{\dbZ} \dbZ_p=\scrF^1\oplus \scrF^0$ such that $\dbG_m$ acts through $\mu_0$ on each $\scrF^i$ via the weight $-i$. 

Let $m\in \dbN^{\ast}$ be defined by the rules: (i) if $(G_1,\scrX_1)$ is of $B_n$ type, then $m:=2n+1$ and (ii) if $(G_1,\scrX_1)$ is of $D_n^{\dbR}$ type, then $m:=2n$. Due to the property 1.4 (iv), $G^{\der}$ is simply connected (cf. also [De2, Table 2.3.8]). Even more, we can also assume that the representation of $G^{\der}_{\dbC}$ on $W\otimes_{\dbQ} \dbC$ is a direct sum of spin representations of direct factors of $G^{\der}_{\dbC}$ which are $\pmb{\text{Spin}}_m$ groups (see loc. cit. and the proof of [De2, Prop. 2.3.10]). This implies that there exists an epimorphism $\theta:G_{\dbZ_p}\twoheadrightarrow J_{\dbZ_p}$ such that the following two properties hold (cf. also [Va3, Subsect. 4.5] for a version of $\theta$ over $\dbR$ in the $D_n^{\dbR}$ type case):

\medskip
{\bf (a)} The kernel of this epimorphism is a closed, flat subgroup scheme of the center of $G_{\dbZ_p}$ (and therefore $\theta$ induces an isomorphism $G_{\dbZ_p}^{\ad}\arrowsim J_{\dbZ_p}^{\ad}$ at the level of adjoint group schemes; if $m=2n+1$, then in fact we have $J_{\dbZ_p}=G_{\dbZ_p}^{\ad}$).
\smallskip
{\bf (b)} The group scheme $J_{\dbZ_p}$ is a finite product $\prod_{i\in I} \Res_{W(k_i)/\dbZ_p} \pmb{\text{SO}}(Q_i,q_i)$ indexed by a finite set $I$ that does not contain $1$, where each $k_i$ is a finite field, $Q_i$ is a free $W(k_i)$-module of rank $m$, and $q_i:Q_i\to W(k_i)$ is a perfect quadratic form. 
\medskip
Let $b_i$ be the perfect bilinear form on $Q_i$ defined by $q_i$; we can view it naturally as a tensor of $Q_i^{\vee}\otimes_{\dbZ_p} Q_i^{\vee}$. Let $(Q,q):=\oplus_{i\in I} (Q_i,q_i)$ be viewed as a free $\dbZ_p$-module $O$ endowed with a perfect quadratic form. Let $\pmb{\text{SO}}_{\dbZ_p}:=\pmb{\text{SO}}(Q,q)$; it is a semisimple group scheme over $\dbZ_p$. Let $\rho:J_{\dbZ_p}\hookrightarrow \pmb{\text{SO}}_{\dbZ_p}$ be the natural faithful representation. Let $(u_{\alpha})_{\alpha\in\scrJ_O}$ be a family of tensors of $\scrT(Q[{1\over p}])$ such that the generic fibre $\pmb{\text{SO}}_{\dbQ_p}$ of $\pmb{\text{SO}}_{\dbZ_p}$ is the subgroup of $\pmb{\text{GL}}_{Q[{1\over p}]}$ that fixes $u_{\alpha}$ for all $\alpha\in\scrJ_O$. We can assume that for each $i\in I$, there exists a subset $\scrJ_{i,O}$ of $\scrJ_O$ that has the following three properties:

\medskip
{\bf (c.i)} each $u_{\alpha}$ with $\alpha\in\scrJ_{i,O}$ is a tensor of $\scrT(Q_i[{1\over p}])$;

\smallskip
{\bf (c.ii)}  the group $\Res_{B(k_i)/\dbQ_p} \pmb{\text{SO}}(Q_i,q_i)_{B(k_i)}$ is the subgroup of $\pmb{\text{GL}}_{Q_i[{1\over p}]}$ that fixes $u_{\alpha}$ for all $\alpha\in\scrJ_{i,O}$;

\smallskip
{\bf (c.iii)}  each element of $W(k_i)\subset \End_{\dbZ_p}(Q_i)$ is of the form $u_{\alpha}$ for some $\alpha\in \scrJ_{i,O}$.  

\medskip
Due to the above description of the $G^{\der}_{\dbC}$-module $W\otimes_{\dbQ} \dbC$, it is easy to see that:

\medskip
{\bf (d)} We can view $Q[{1\over p}]$ as a $G_{\dbQ_p}$-submodule of $\End_{\dbQ_p}(W^{\vee}\otimes_{\dbQ} \dbQ_p)\subseteq\scrT(W^{\vee}\otimes_{\dbQ} \dbQ_p)$. 

\medskip
In other words, the standard $m$-dimensional representation of an $\pmb{{\grs}{\gro}}_m$ Lie algebra over $\dbC$ (or over $\dbQ_p$) is a direct summand of the representation of $\pmb{{\grs}{\gro}}_m$ which is the tensor product of the spin representation and of the dual of the spin representation. 

Due to the property (d), we can view each tensor $u_{\alpha}$ as a linear combination with coefficients in $\dbQ_p$ of the family $(v_{\alpha})_{\alpha\in\scrJ}$ of tensors of $\scrT(W^{\vee})\subset\scrT(W^{\vee}\otimes_{\dbQ} \dbQ_p)$ introduced in Subsubsection 2.2.1. Let $\mu_i:\dbG_m\to \Res_{W(k_i)/\dbZ_p} \pmb{\text{SO}}(Q_i,q_i)$ be the cocharacter induced naturally by $\mu_0$ via the epimorphism $\theta:G_{\dbZ_p}\twoheadrightarrow J_{\dbZ_p}$. We have an extra property:

\medskip
{\bf (e)} If $p>2$, then for each $i\in I$ the family of tensors of $\scrT(Q_i)$ formed by $b_i$ and $W(k_i)$, is strongly $\dbZ_p$-very well position for the closed subgroup scheme $\Res_{W(k_i)/\dbZ_p} \pmb{\text{SO}}(Q_i,q_i)$ of $\pmb{\text{GL}}_{O_i}$ in the sense of [Va1, Def. 4.3.4 and Rm. 4.3.7 1)].

\medskip
In other words for $p>2$, if $C$ is a flat, reduced $\dbZ_p$-algebra and if $\tilde Q_i$ is a free $C$-module such that (i) we have $\tilde Q_i[{1\over p}]=Q_i\otimes_{\dbZ_p} C[{1\over p}]$, (ii) $b_i$ induces a perfect bilinear form on $\tilde Q_i$, and (iii) $W(k_i)\otimes_{\dbZ_p} C$ is a $C$-subalgebra of $\End_C(\tilde Q_i)$, then the schematic closure of $\Res_{W(k_i)/\dbZ_p} \pmb{\text{SO}}(Q_i,q_i)\times_{\Spec\dbZ_p} \Spec C[{1\over p}]$ in $\pmb{GL}_{\tilde Q_i}$ is a reductive group scheme. To check this, we can assume that the $C$-algebra $W(k_i)\otimes_{\dbZ_p} C$ is isomorphic to $C^{[k_i:\dbF_p]}$ and this case is standard. 

The dimension of $\scrX_1$ as a complex manifold (i.e., the dimension of the Shimura variety $\Sh(G_1,\scrX_1)$) is $m-2$. We will use the notations of Subsubsection 2.2.1. By enlarging the algebraically closed field $k$, we can assume that there exists a finite, discrete valuation ring extension $V$ of $W(k)$ such that the normalization $\scrP$ of $\scrN_V$ is regular in characteristic $0$ and in codimension $1$ (cf. [PY, Appendix] and the fact that $\scrN$ is a pro-\'etale cover of a quasi-projective scheme over $\dbZ_{(p)}$). Let $K$ be the field of fractions of $V$. Let $e$ be the index of ramification of $V$ and let $\pi$ be a uniformizer of $V$. The minimal polynomial $h_e\in W(k)[X]$ of $\pi$ over $W(k)$ is an Eisenstein polynomial.

For a perfect field $k_1$ that contains $k$, let $\grS_1:=W(k_1)[[x]]$ and let $R_{e,1}$ be the $\grS_1$-subalgebra of $B(k_1)[[x]]$ formed by formal power series $\sum_{n\ge 0} a_nx^n$ such that the sequence $([{n\over e}]!a_n)_{n\in\dbN\cap\{0\}}$ is formed by elements of $W(k_1)$ and converges to $0$ in the $p$-adic topology of $W(k_1)$; here $x$ is an independent variable. We have a $W(k_1)$-epimorphism $m_V:R_{e,1}\twoheadrightarrow V\otimes_{W(k)} W(k_1)$ that maps $x$ to $\pi\otimes 1$. Let $\Phi_{k_1}$ be the Frobenius lift of $\grS_1$ or $R_{e,1}$ that is compatible with $\sigma_{k_1}$ and that maps $x$ to $x^p$.

Let $\scrO$ be the open subscheme of $\scrP_k$ which is the ordinary locus: a point $y\in\scrP(k)$ belongs to $\scrO(k)$ if and only if the abelian variety $y^*(\scrA_{\scrP})$ is ordinary.

\medskip\noindent
{\bf 3.6.1. Proposition.} {\it The ordinary locus $\scrO$ is Zariski dense in $\scrP_k$.}

\medskip
\proof
It suffices to show that if $Y=\Spec R$ is an arbitrary non-empty, affine, connected, regular $k$-subscheme of $\scrP_k$ of dimension $n$, then the abelian scheme $\scrA_Y$ is generically ordinary. To check this we can perform the following two replacement operations:

\medskip
{\bf (o1)} Replace $Y$ by an affine, open, Zariski dense subscheme of it.

\smallskip
{\bf (o2)} Replace $(\scrP,Y)$ by $(\scrP^\prime,Y\times_{\scr P} \scrP^\prime)$, where $\scrP^\prime$ is an affine, quasi-finite scheme over $\scrP$ which is regular and formally smooth over $V$ and whose special fibre $\scrP_k^\prime$ dominates $\scrP_k$. 
\medskip

By performing the operation (o1) we can assume that there exists a smooth, affine $W(k)$-scheme $\Spec \scrR$ whose reduction modulo $p$ is $Y$ and for which there exists an ind-\'etale $W(k)$-homomorphism $W(k)[x_1,\ldots,x_{m-2}]\to\scrR$. Let $\scrR^\wedge$ be the $p$-adic completion of $\scrR$ and let $\Phi_{\scrR}$ be the Frobenius lift of it which is compatible with $\sigma$ and which takes $x_i$ to $x_i^p$ for all $i\in\{1,\ldots,m-2\}$. We can assume that $\scrZ:=\Spec V\otimes_{W(k)} \scrR^\wedge$ is the affine scheme defined naturally by the $p$-adic completion of an affine, open subscheme of $\scrP^\prime$ which lifts $Y$; thus we can speak about the abelian scheme $\scrA_\scrZ$.

Let $\scrC$ be the $F$-isocrystal over $R$ of the $p$-divisible group of $\scrA_Y$. It is defined naturally by a projective $\scrR^\wedge[{1\over p}]$-module $\scrV$ equipped with a $\Phi_{\scrR}$-linear endomorphism $\phi_{\scrV}:\scrV\to \scrV$. The $K\otimes_{B(k)} \scrR^\wedge$-module $K\otimes_{B(k)} \scrV$ is equipped with a direct summand $\scrF$ defined by the Hodge filtration of  the $p$-divisible group of $\scrA_\scrZ$. In what follows we will express this property by saying that $\scrC$ gets a filtration after tensorization over $B(k)$ with $K$. Accordingly, below all $F$-isocrystals over (some completion of) $R$ will get similar filtrations after tensorization over $B(k)$ with $K$ and all morphisms of $F$-isocrystals over (some completion of) $R$ will be compatible with the corresponding filtrations one gets after tensorization over $B(k)$ with $K$. Moreover all the Kodaira--Spencer maps of $F$-isocrystals over (some completion of) $R$, will be computed after tensorization over $B(k)$ with $K$ and will be with respect to the mentioned filtrations one gets  after tensorization over $B(k)$ with $K$.

Let $\kappa$ be the field of fractions of $R$. Let $O$ be the $p$-adic completion of the local ring of $\scrR$ (or of $\scrR^\wedge$) whose residue field is $\kappa$. It is a discrete valuation ring of mixed characteristic $(0,p)$ and index of ramification $1$. Let $V_1:=V\otimes_{W(k)} O$ and let $K_1$ be the field of fractions of $V_1$. We denote by $z_1$ the composite morphism $\Spec  V_1\to \scrZ\to\scrP^\prime\to \scrN$. Let $\rho_1:\Gal(K_1)\to G_{\dbZ_p}(\dbZ_p)$ be the homomorphism associated naturally to the $p$-adic Galois representation of the abelian scheme $A_1=z_1^*(\scrA)$ over $V_1$; this makes sense based on [Va15, Lem. 2.3.4 (a)]. By composing $\rho_1$ with $\theta(\dbZ_p):G_{\dbZ_p}(\dbZ_p)\to J_{\dbZ_p}(\dbZ_p)$, we get a homomorphism 
$$\rho_1^{\pmb{\text{SO}}}:\Gal(K_1)\to J_{\dbZ_p}(\dbZ_p)=\prod_{i\in I} \pmb{\text{SO}}(O_i,q_i)(W(k_i)).$$ 
\indent
Let $\scrI=\oplus_{i\in I} \scrI_i$ be the $F$-subisocrystal of $\scrT(\scrC):=\oplus_{s,t\in\dbN} \scrC^{\otimes s}\otimes \scrC^{{\vee}\otimes t}$ (equivalently, of $End(\scrC)$) that corresponds naturally to the $G_{\dbQ_p}$-submodule $Q[{1\over p}]=\oplus_{i\in I} Q_i[{1\over p}]$ of $\scrT(W^{\vee}\otimes_{\dbQ} \dbQ_p)$ (equivalently, of $\End_{\dbQ_p}(W^{\vee}\otimes_{\dbQ} \dbQ_p)$). In what follows, by a {\it latticed $F$-isocrystal} we mean a Tate-twist $\ddag(s)$ of an $F$-crystal $\ddag$; here $s\in\dbZ$. 

By performing the operations (o1) and (o2),  we can assume that each $\scrI_i$ is the $F$-isocrystal of a latticed $F$-isocrystal $\scrQ_i$ over $R$ which generically (i.e., over $\kappa$) corresponds naturally to the homomorphism $\rho_{1,i}^{\pmb{\text{SO}}}:\Gal(K_1)\to \pmb{\text{SO}}(O_i,q_i)(W(k_i))$ defined naturally by $\rho_1^{\pmb{\text{SO}}}$. Even more, by performing the operations (o1) and (o2) we can assume that for each element $i\in I$ the following three properties hold:

\medskip
{\bf (a)} The latticed $F$-isocrystal $\scrQ_i$ over $R$ is equipped naturally with a perfect quadratic form $\scrK_i$ (that corresponds to $q_i$ via Fontaine comparison theory).

\smallskip
{\bf (b)} We can identify $(\scrQ_i,\scrK_i)$ with
$$\scrE_i:=(Q_i\otimes_{\dbZ_p} \scrR^\wedge,g_{i,\scrR}(\mu_i({1\over p})\otimes\Phi_{\scrR}),q_i,\nabla_i)$$
for some element $g_{i,\scrR}\in \Res_{W(k_i)/\dbZ_p} \pmb{\text{SO}}(Q_i,q_i)(\scrR^\wedge)$ and some integrable, nilpotent modulo $p$ connection $\nabla_i$ on $Q_i\otimes_{\dbZ_p} \scrR^\wedge$. Under this identification, each tensor $u_{\alpha}\in\scrT(Q_i\otimes_{\dbZ_p} \scrR^\wedge[{1\over p}])$ generates the $F$-isocrystal over $R$ that corresponds naturally to the tensor $u_{\alpha}\in\scrT(Q_i[{1\over p}])$.

\smallskip
{\bf (c)} The resulting filtration of $K\otimes_{B(k)} (Q_i\otimes_{\dbZ_p}\scrR^\wedge[{1\over p}])$ (induced via (b) from the one with which $\scrQ_i$ is equipped after tensorization over $B(k)$ with $K$) is such that it induces a filtration of $V\otimes_{W(k)} (Q_i\otimes_{\dbZ_p} \scrR^\wedge)$ defined by a cocharacter 
$$\mu_{i,V}:\dbG_m\to [\Res_{W(k_i)/\dbZ_p}\pmb{\text{SO}}(Q_i,q_i)]_{V\otimes_{W(k)} \scrR^\wedge}$$ 
that lifts the extension to $V$ of the cocharacter $\mu_i:\dbG_m\to\Res_{W(k_i)/\dbZ_p}\pmb{\text{SO}}(Q_i,q_i)$.

\medskip
To check that all these three properties hold it suffices to check that they hold over $\bar V_1:=V/pV\otimes_k k_1$, where $k_1$ is the perfection of $\kappa$. Here are two ways to argue that all these three properties hold over $\bar V_1$ (the first one works only for $p\ge 5$). 

\medskip
{\bf (i)} For $p\ge 5$, the fact that we can assume that properties (a) to (c) hold over $\bar V_1$ (resp. directly over $V\otimes_{W(k)} \scrR^\wedge$) follows from [Fa, Thm. 5 iii)] applied over the discrete valuation ring $V\otimes_{W(k)} W(k_1)$ that dominates $O_1$ (resp. from [Fa, Thm. 5*]). Based on the property 3.6 (e), the arguments are the same as the ones of [Va1, Subsect. 5.2].

\smallskip
{\bf (ii)} For $p\ge 2$, the fact that we can assume that properties (a) and (b) hold (resp. that the property (c) holds) over $\bar V_1$ follows from [Ki3, Prop. (1.3.4) and Cor. (1.3.5)] (resp. from [Ki3, Prop. (1.1.5) and Lem. (1.4.5)]; to be compared with [Va1, Subsubsects. 5.3.1 and 5.3.2] and [Va11, Lem. 5.2.6]). To detailed on this, let $(\grQ_i,p^{-1}\varphi_i)$ be the contravariant {\it Breuil--Kisin module} associated to the Galois representation $\rho_{1,i,k_1}^{\pmb{\text{SO}}}:\Gal(K\otimes_{B(k)} B(k_1))\to \pmb{\text{SO}}(O_i,q_i)(W(k_i))$ induced naturally by $\rho_{1,i}^{\pmb{\text{SO}}}$. We recall that $\grQ_i$ is a free $\grS_1$-module of the same rank as $Q_i$ and that $\varphi_i:\grQ_i\otimes_{\grS_1} {}_{\Phi_{k_1}} \grS_1\to\grQ_i$ is a $\grS_1$-linear map whose cokernel is annihilated by $h^2_e$. To each $u_{\alpha}$ with $\alpha\in\scrJ_{i,O}$ corresponds naturally a tensor $t_{\alpha}\in\scrT(\grQ_i[{1\over p}])$ (cf. [Ki1]). As $\Res_{W(k_i)/\dbZ_p} \pmb{\text{SO}}(O_i,q_i)$ is a semisimple (and thus reductive) group scheme over $\dbZ_p$, from [Ki3, Prop. (1.3.4) and Cor. (1.3.5)] we get that we have an isomorphism $(\grQ_i,(t_{\alpha})_{\alpha\in\scrJ_{i,O}})\arrowsim (Q_i\otimes_{\dbZ_p} \grS_1,(u_{\alpha})_{\alpha\in\scrJ_{i,O}})$ and thus (cf.  property 3.6 (c.iii)) $\grQ_i$ is a $W(k_i)$-module equipped with a perfect quadratic form. The canonical way of passing from Breuil-Kisin modules to $F$-crystals (see [Ki1--3]) implies that $(\grQ_i,p^{-1}\varphi_i,(u_{\alpha})_{\alpha\in\scrJ_{i,O}})\otimes_{\grS_1} {}_{\Phi_{k_1}} R_{e,1}$ is the latticed $F$-isocrystal endowed with tensors over $\bar V_1$ which induces naturally the searched for latticed $F$-isocrystal $\scrQ_{i,k_1}$ endowed with tensors over $k_1$. Thus the fact that properties (a) and (b) hold (resp. that the property (c) holds) over $\bar V_1$ follows from the last two sentences (resp. from [Ki3, Prop. (1.1.5) and Lem. (1.4.5)]).

\medskip
Let $\scrR_0$ be the completion of $\scrR^\wedge$ at some $k$-valued point of it. Let $\Phi_0$ be the Frobenius lift of $\scrR_0$ induced naturally by $\Phi_{\scrR}$. We can assume that $\scrR_0=W(k)[[x_1,\ldots,x_{m-2}]]$ and that $\Phi_0$ takes $x_j$ to $x_j^p$ for all $j\in\{1,\ldots,n\}$. Let $g_{i,0}\in\Res_{W(k_i)/\dbZ_p}\pmb{\text{SO}}(Q_i,q_i)(W(k))$ be the reduction modulo $(x_1,\ldots,x_{m-2})$ of $g_{i,\scrR}$. Let 
$$g_0:=(g_{i,0})_{i\in I}\in J_{\dbZ_p}(W(k))=\prod_{i\in I}\Res_{W(k_i)/\dbZ_p}\pmb{\text{SO}}(Q_i,q_i)(W(k)).$$
\indent
Let $U_0$ be the maximal, unipotent, smooth, closed subgroup of either $G_{\dbZ_p}$ or $J_{\dbZ_p}$ with the property that $\mu_0$ acts identically on $\Lie(U_0)$. We can identify $\scrR_0$ with the local ring of the completion of $U_0$ along its identity section. Thus we have a universal element $u_0\in U_0(\scrR_0)$. Let 
$$\scrE_0:=(Q\otimes_{\dbZ_p} \scrR_0,u_0(g_0\mu_0({1\over p})\otimes\Phi_0),q,\nabla_0),$$
where $\nabla_0:Q\otimes_{\dbZ_p} \scrR_0\to Q\otimes_{\dbZ_p} (\oplus_{j=1}^{m-2} dx_i\scrR_0)$ is an integrable, nilpotent modulo $p$ connection on $Q\otimes_{\dbZ_p} \scrR_0$ which satisfies the following  identity
$$\nabla_0\circ (u_0(g_0\mu_0({1\over p})\otimes\Phi_0))=(u_0(g_0\mu_0({1\over p})\otimes\Phi_0)\otimes d\Phi_{\scrR_0})\circ\nabla_0.\leqno (10)$$
It is well known that there is at most one connection $\nabla_0$ satisfying the identity (10) (the argument for this is the same as the one of [Va11, Subsect. 5.2.1])  and below we will argue that the connection $\nabla_0$ exists.

We claim that the pull-back $(Q\otimes_{\dbZ_p} \scrR_0,u_1(g_0\mu_0({1\over p})\otimes\Phi_{R_0}),q,\nabla_1)$ of $\oplus_{i\in I} \scrE_i$ to $\scrR_0/p\scrR_0$ is the pull-back of $\scrE_0$ via a unique morphism $q_0:\Spec \scrR_0/p\scrR_0\to\Spec \scrR_0/p\scrR_0$; here $u_1\in J_{\dbZ_p}(\scrR_0)$ is congruent to the identity element modulo the ideal $(x_1,\ldots,x_{m-2})$ of $\scrR_0$. If $p>2$, then there exists an element $g\in G_{\dbZ_p}(W(k))$ whose image in $G_{\dbZ_p}^{\ad}(W(k))$ is $g_0$. If $p=2$, then the element $g$ exists provided we replace $g_0$ by the image in $J_{\dbZ_p}(W(k))$ of a suitable element $hg\phi(h^{-1})$, where $h\in G_{\dbZ_p}(W(k))$ normalizes $\scrF^1/p\scrF^1$ (to be compared with [Va10, Fact 2.6.3]). A similar argument shows that we can assume that there exists an element $g_1\in G_{\dbZ_p}(\scrR_0)$ that maps to $u_1$. 

As the connection $\nabla_0$ (resp. $\nabla_1$) is uniquely determined by $u_0$ (resp. $u_1$) via the identity (10) (resp. the identity $\nabla_1\circ (u_1(g_0\mu_0({1\over p})\otimes\Phi_{R_0}))=(u_1(g_0\mu_0({1\over p})\otimes\Phi_{R_0})\otimes d\Phi_{\scrR_0})\circ\nabla_1$), to prove the claim it suffices to show that $(L^{\vee}_p\otimes_{\dbZ_p} \scrR_0,\scrF^1\otimes_{\dbZ_p} \scrR_0,g_1(g\mu_0({1\over p})\otimes\Phi_{R_0}),(v_{\alpha})_{\alpha\in\scrJ},\nabla_1^\prime)$ is the pull-back of $(L^{\vee}_p\otimes_{\dbZ_p} \scrR_0,\scrF^1\otimes_{\dbZ_p} \scrR_0,u_0(g\mu_0({1\over p})\otimes\Phi_{R_0}),(v_{\alpha})_{\alpha\in\scrJ},\nabla_0^\prime)$ via a unique morphism $\Spec \scrR_0/p\scrR_0\to\Spec \scrR_0/p\scrR_0$  which at the level of rings maps the ideal $(x_1,\ldots,x_{m-2})$ into itself; here the connections $\nabla_1^\prime$ and $\nabla_0^\prime$ are uniquely determined by the identities $\nabla_1^\prime\circ (g_1(g\mu_0({1\over p})\otimes\Phi_{R_0}))=(g_1(g\mu_0({1\over p})\otimes\Phi_{R_0})\otimes d\Phi_{\scrR_0})\circ\nabla_1^\prime$ and $\nabla_0^\prime\circ (u_0(g\mu_0({1\over p})\otimes\Phi_{R_0}))=(u_0(g\mu_0({1\over p})\otimes\Phi_{R_0})\otimes d\Phi_{\scrR_0})\circ\nabla_0^\prime$ (respectively) --cf. [Va15, Appendix, Subsect. B6]-- and induce respectively the connections $\nabla_1$ and $\nabla_0$ on $Q\otimes_{\dbZ_p} \scrR_0$ (this proves as well the existence of $\nabla_0$). But this pull-back propriety is a particular case of [Va15, Appendix, Thm. B6.4]. Thus the claim holds. 

The Kodaira--Spencer map of the connection $\nabla_0^\prime$ (and thus also of $\nabla_0$) is injective and its image is canonically identified with $\Lie(U_0)\otimes_{\dbZ_p} \scrR_0$, cf. [Va15, Appendix, property B6.3 (iii)]. Due to this and the property (c), it is easy to see that the morphism $q_0$ lifts to a morphism $q:\Spec V\otimes_{W(k)} \scrR_0\to \Spec \scrR_0$ such that the following property holds:

\medskip
{\bf (d)} for each $i\in I$  the filtration of $V\otimes_{W(k)} (Q_i\otimes_{\dbZ_p} \scrR^\wedge)$ induced by the cocharacter $\mu_V$ of the property (c) is the pull-back via $q$ of the filtration of $V\otimes_{W(k)} (Q_i\otimes_{\dbZ_p} \scrR^\wedge)$ induced by $\mu_{i,V\otimes_{W(k)} (Q_i\otimes_{\dbZ_p} \scrR^\wedge)}$.

\medskip
We show that the assumption that the morphism $q_0:\Spec \scrR_0/p\scrR_0\to\Spec \scrR_0/p\scrR_0$ is not dominant leads to a contradiction. From this assumption, from the fact that the connection $\nabla_0$ is versal, and from the property (d)  we get that there exists an element $i_0\in I$ such that the Kodaira--Spencer map of $\nabla_{i_0}$ has an image whose rank is less than the relative dimension of $U_{i_0,0}:=U_0\cap \Res_{W(k_{i_0})/\dbZ_p}\pmb{\text{SO}}(Q_{i_0},q_{i_0})$. But the Kodaira--Spencer map of the pull-back $\scrC_{\scrR_0/p\scrR_0}$ of $\scrC$ to $\scrR_0/p\scrR_0$ can be identified with the tensorization with $K$ over $B(k)$ of the direct sum of the Kodaira--Spencer maps of $\nabla_i$'s with $i\in I$ (this is so as the adjoint representation of $G_{\dbZ_p}$ is the composite of the homomorphism $\theta:G_{\dbZ_p}\twoheadrightarrow J_{\dbZ_p}$ with the adjoint representation of $J_{\dbZ_p}$). From the last two sentences we get that the Kodaira--Spencer map of $\scrC_{\scrR_0/p\scrR_0}$ is not injective. As $Y$ is a quasi-finite over a pro-\'etale cover of a scheme which is finite over $\scrA_{d,1,N,k}$, the Kodaira--Spencer map of $\scrC_{\scrR_0/p\scrR_0}$ is injective. Contradiction. Therefore the morphism $q_0$ is dominant. 

We claim that $\scrE_0$ is generically ordinary. But this is a direct consequence of the following two statements (the first one being trivial):

\medskip
{\bf (e)} the latticed $F$-isocrystal $(Q,\mu_0({1\over p}))$ over $\dbF_p$ is ordinary and there exists a largest Zariski dense, open subscheme $\scrU$ of the special fibre $J_{\dbF_p}$ of $J_{\dbZ_p}$ such that for each geometric point $\bar h: \Spec K\to \scrU$ and for every element $h\in J_{\dbZ_p}(W(K))$ that lifts $\bar h$, the latticed $F$-isocrystal $(Q\otimes_{\dbZ_p} W(K),h(\mu_0({1\over p})\otimes\sigma_K))$ over $K$ is ordinary;

\smallskip
{\bf (f)} the closed subscheme $U_{0,\dbF_p}$ of $J_{\dbF_p}$ has a non-empty intersection with $\scrU g_0^{-1}$.

\medskip
To check the property (f), let $\scrW$ be the parabolic subgroup of $J_{\dbF_p}$ which is the image of the parabolic subgroup of the special fibre $G_{\dbF_p}$ of $G_{\dbZ_p}$ which is the normalizer of $\scrF^1/p\scrF^1$ in $G_{\dbF_p}$. If $l=l_0l_1\in J_{\dbZ_p}(W(K))$ lifts a $K$-valued point of $\scrW$ and $l_0$ is fixed by $\mu_0$ and $\dbG_m$ acts through $\mu_0$ on the $\dbZ_p$-span of $l_1-1_{Q}$ via the weight $-1$, then for each element $h\in J_{\dbZ_p}(W(K))$, the latticed $F$-isocrystals $(Q\otimes_{\dbZ_p} W(K),lh(\mu_0({1\over p})\otimes\sigma_K)l^{-1})=(Q\otimes_{\dbZ_p} W(K),lh\sigma_K(\tilde l_0l_1^p)^{-1}(\mu_0({1\over p})\otimes\sigma_K))$ and $(Q\otimes_{\dbZ_p} W(K),h(\mu_0({1\over p})\otimes\sigma_K))$ are isomorphic. Thus to check (f) it suffices to show that the morphism $s_0:\scrW\times_{\dbF_p} U_{0,\dbF_p}\to G_{\dbF_p}$ which at the level of $K$-valued points maps the pair $(\bar l=\bar l_0\bar l_1,\bar u)$ to $\bar l\bar u \bar g_0\sigma_K(\bar l_0)^{-1}$ (with $\bar g_0$ as $g_0$ modulo $p$), has an open image. But this is a direct consequence of the fact that the tangent map of $s_0$ at the identity element of $\scrW\times_{\dbF_p} U_{0,\dbF_p}\to G_{\dbF_p}$ is a bijection (we note that the product morphism $\scrW\times_{\dbF_p} U_{0,\dbF_p}\to G_{\dbF_p}$ is an open embedding). Thus (f) holds. 

As $q_0$ is dominant and as $\scrE_0$ is generically ordinary, by performing the operation (o1), we can assume that $\scrI$ is an ordinary $F$-isocrystal and therefore that $\scrC$ is ordinary. This implies that $\scrA_Y$ is generically ordinary. Thus the proposition holds.\endproof 

\medskip\noindent
{\bf 3.6.2. Conclusion.} From [Va15, Thm. 1.7 (a) and (b)] and Proposition 3.6.1 we get that $\scrN^{\s}=\scrN$. Thus, regardless of what $p$ is, the property 1.5 (i) holds in this last Case 6. If $p>2$, then as in the last part of the Case 3 (Subsection 3.3) we argue that Theorem 1.5 (ii) and (iii) hold in this last Case 6.

\bigskip\noindent
{\bf 3.6.3. Remark.} If $p\ge 5$, then we can choose $(\scrQ_i,\scrK_i)$ to be defined globally on $\scrP_k=\scrN_k$ and in fact to correspond canonically to $(Q_i,q_i)$ via Fontaine comparison theory (cf. [Fa, Thm. 5*]). 

\bigskip\noindent
{\bf 3.7. Proposition.} {\it Suppose that $(G_1,\scrX_1,H_1)$ is a Shimura triple of abelian type with respect to $p$ such that the Shimura pair $(G_1,\scrX_1)$ is simple, adjoint. Then there exists a commutative diagram of Shimura triples of abelian type 
$$ 
\spreadmatrixlines{1\jot}
\CD
(G_4,\scrX_4,H_4) @>{f_3}>> (G_3,\scrX_3,H_3)\\
@V{\pi_1}VV @VV{\pi_3}V\\
(G_1,\scrX_1,H_1) @>{f_1}>> (G_2,\scrX_2,H_2)
\endCD
$$
such that the following four properties hold:

\medskip
{\bf (i)} the Shimura pair $(G_2,\scrX_2)$ is adjoint and unitary;

\smallskip
{\bf (ii)} both horizontal maps $f_1$ and $f_3$ are injective;

\smallskip
{\bf (iii)} both vertical maps $\pi_1$ and $\pi_3$ induce isomorphisms at the level of adjoint Shimura triples (i.e., they induce naturally isomorphisms $(G_4^{\ad},\scrX_4^{\ad},H_4^{\ad})\arrowsim (G_1,\scrX_1,H_1)$ and $\break(G_3^{\ad},\scrX_3^{\ad},H_3^{\ad})\arrowsim (G_2^{\ad},\scrX_2^{\ad},H_2^{\ad})$);

\smallskip
{\bf (iv)} the derived group $G_4^{\der}$ is the maximal one allowed by the abelian type.}

\medskip
\proof
We can assume that $(G_1,\scrX_1)$ is not unitary. This proposition is only a $\dbZ_{(p)}$ version of the results of Satake on embeddings between hermitian symmetric domains of classical Lie type (see [Sa1,2]). If $(G_1,\scrX_1)$ is not (resp. is) of $D_n^{\dbR}$ type, then a $\dbQ$--version of this proposition is presented in [Va8, Subsects. 4.2 to 4.8] (resp. in loc. cit. and [Va8, Rm. 4.8.2 (c)]). The passage from $\dbQ$-versions to $\dbZ_{(p)}$-versions is standard and the easy details are left to the reader (they are also presented in [Va16, Subsect. 5.1 (a) to (g)]).\endproof 

\bigskip\noindent
{\bf 3.8. End of the proof of Theorem 1.5.}
We know that Theorem 1.5 (i) holds, cf. Subsections 3.1 to 3.5 and Subsubsection 3.6.2. We now check that Theorem 1.5 (ii) and (iii) hold. We can assume that $(G_1,\scrX_1)$ is not unitary, cf. Subsection 3.1. If $p>2$, then Theorem 1.5 (ii) and (iii) follow from Subsections 3.2, 3.3, and 3.5 and from Subsubsection 3.6.2. We will use the notations of Propositions 1.4 and  3.7 in order to show that Theorem 1.5 (ii) and (iii) holds even if $p=2$. For the sake of uniformity, we will continue to work with an arbitrary prime $p$. The central isogenies $G^{\der}\to G_1$ and $G_4^{\der}\to G_1$ can be identified, cf. properties 1.4 (iv) and 3.7 (iv). 

Let $\scrN_2$ and $\scrN_3$ be the integral canonical models of $(G_2,\scrX_2,H_2)$ and $(G_3,\scrX_3,H_3)$ (respectively), cf. [Va7, Thm. 1.3]. They are quasi-projective, cf. loc. cit.  As the integral canonical model $\scrN$ of $(G,\scrX,H)$ exists (cf. Theorem 1.5 (i)) and is quasi-projective, the integral canonical model $\scrN_4$ of $(G_4,\scrX_4,H_4)$ exists as well (cf. [Va7, Prop. 2.4.3 (a)]). From Proposition 2.4 we get that $\scrN_4$ is the normalization of the schematic closure of $\Sh(G_4,\scrX_4)/H_4$ in $\scrN_{3,E(G_4,\scrX_4)_{(p)}}$  (via the functorial closed embedding morphism $f_3:\Sh(G_4,\scrX_4)/H_4\hookrightarrow \Sh(G_3,\scrX_3)_{E(G_4,\scrX_4)}/H_3$).

Let $\scrN_1$ be the normalization of the schematic closure of $\Sh(G_1,\scrX_1)/H_1$ in $\scrN_{2,E(G_1,\scrX_1)_{(p)}}$ (via the functorial closed embedding morphism $f_1:\Sh(G_1,\scrX_1)/H_1\hookrightarrow \Sh(G_2,\scrX_2)_{E(G_1,\scrX_1)}/H_2$). It is a normal integral model of $(G_1,\scrX_1,H_1)$ which has the extension property, cf. [Va1, Prop. 3.4.1]. As $\scrN_2$ is quasi-projective, $\scrN_1$ is also quasi-projective. Due to the extension properties enjoyed by $\scrN_1$ to $\scrN_4$ and the fact that $\scrN_2$ to $\scrN_4$ are healthy regular schemes (being regular and formally smooth over $\dbZ_{(p)}$), we have a commutative diagram
$$ 
\spreadmatrixlines{1\jot}
\CD
\scrN_4 @>{f_{3,p}}>> \scrN_{3,E(G_4,\scrX_4)_{(p)}}\\
@V{\pi_{1,p}}VV @VV{\pi_{3,p}}V\\
\scrN_{1,E(G_4,\scrX_4)_{(p)}} @>{f_{1,p}}>> \scrN_{2,E(G_4,\scrX_4)_{(p)}}
\endCD
$$ 
of normal $E(G_4,\scrX_4)_{(p)}$-schemes. The following properties hold:

\medskip
{\bf (v)} the morphism $\pi_{3,p}$ is a pro-\'etale cover of an open closed subscheme of $\scrN_{2,E(G_4,\scrX_4)_{(p)}}$;

\smallskip
{\bf (vi)} the morphisms $f_{3,p}$, $\pi_{1,p}$, and $f_{1,p}$ are pro-finite;

\smallskip
{\bf (vii)} the generic fibre of $\pi_{1,p}$ is a pro-\'etale cover of an open closed subscheme of $\scrN_{1,E(G_4,\scrX_4)}$.

\medskip
The property (v) is implied by [Va7, Thm. 5.1 (b)]. Due to properties (vi) and (vii), the image of $\pi_{1,p}$ is an open closed subscheme of $\scrN_{1,E(G_4,\scrX_4)_{(p)}}$. Due to the property (v) we easily get that $\pi_{1,p}$ is in fact a pro-\'etale cover of its image; thus this image is a regular, formally smooth scheme over $E(G_4,\scrX_4)_{(p)}$ and thus also over $E(G_1,\scrX_1)_{(p)}$. The connected components of $\scrN_{1,E(G_4,\scrX_4)_{(p)}}$ are permuted transitively by $G_1(\dbA_f^{(p)})$, cf. [Va1, Lem. 3.3.2]. From the last two sentences we get that $\scrN_{1,E(G_4,\scrX_4)_{(p)}}$ is a regular, formally smooth scheme over $E(G_4,\scrX_4)_{(p)}$. This implies that $\scrN_1$ is a regular, formally smooth scheme over $E(G_1,\scrX_1)_{(p)}$. From this and [Va1, Cor. 3.4.4] we get that $\scrN_1$ is the integral canonical model of $(G_1,\scrX_1,H_1)$. Thus Theorem 1.5 (ii) holds. Moreover we have a functorial morphism $\scrN\to\scrN_1$ of $E(G_1,\scrX_1)_{(p)}$-schemes.

Let $\dbZ_{(p)}^{\text{un}}$ be the maximal $\dbZ_{(p)}$-subalgebra of $\overline{\dbQ}$ such that $\Spec \dbZ_{(p)}^{\text{un}}$ is a pro-\'etale cover of $\Spec \dbZ_{(p)}$. Both $E(G_4,\scrX_4)_{(p)}$ and $E(G,\scrX)_{(p)}$ are $\dbZ_{(p)}$-subalgebras of $\dbZ_{(p)}^{\text{un}}$ and moreover the connected components of $\scrN_{\dbZ_{(p)}^{\text{un}}}$ and $\scrN_{4,\dbZ_{(p)}^{\text{un}}}$ can be canonically identified, cf. [Va7, Prop. 2.4.3 (c)]. As the connected components of $\scrN_{4,\dbZ_{(p)}^{\text{un}}}$ are pro-\'etale covers of certain connected components of $\scrN_{1,\dbZ_{(p)}^{\text{un}}}$, we get that there exists a connected component of $\scrN$ which is a pro-\'etale cover of an open closed subscheme of $\scrN_1$. As the connected components of $\scrN$ are permuted transitively by $G(\dbA_f^{(p)})$ (cf. [Va1, Lem. 3.3.2]), we get that $\scrN$ itself is a pro-\'etale cover of an open closed subscheme of $\scrN_1$. Thus Theorem 1.5 (iii) holds as well. This ends the proof of Theorem 1.5.\endproof

\medskip
We have the following corollary to Theorem 1.5. 

\bigskip\noindent
{\bf 3.9. Corollary.} {\it Let $(G_1,\scrX_1,H_1)$ and $\scrN_1$ be as in Theorem 1.5. Let $f_2:(G_2,\scrX_2,H_2)\to (G_1,\scrX_1,H_1)$ be a map of Shimura triples of abelian type that induces an isomorphism $(G_2^{\ad},\scrX_2^{\ad},H_2^{\ad})\arrowsim (G_1,\scrX_1,H_1)$. Then the normalization $\scrN_2$ of $\scrN_1$ in the ring of fractions of $\Sh(G_2,\scrX_2)/H_2$ is a pro-\'etale cover of an open closed subscheme of $\scrN_1$. Moreover, $\scrN_2$ together with the natural action of $G_2(\dbA_f^{(p)})$ on it, is the integral canonical model of $(G_2,\scrX_2,H_2)$ and it is quasi-projective. If  the integral canonical model $\scrN$ of Theorem 1.5 is projective, then $\scrN_2$ is projective too.}

\medskip
\proof
We will use the notations of Theorem 1.5. We consider the fibre product of $f:(G,\scrX,H)\to (G_1,\scrX_1,H_1)$ and $f_2$ (cf. [Va1, Subsect. 2.4 and Rm. 3.2.7 3)])

$$
\spreadmatrixlines{1\jot}
\CD
(G_3,\scrX_3,H_3) @>{\pi}>> (G,\scrX,H)\\
@V{\pi_1}VV @VV{f}V\\
(G_2,\scrX_2,H_2) @>{f_2}>> (G_1,\scrX_1,H_1).
\endCD
$$
Due to the property 1.4 (iv), we have $G_3^{\der}=G^{\der}$. Thus by applying [Va7, Prop. 2.4.3 (a) and (b)] to $(G_3,\scrX_3,H_3)$ and $(G,\scrX,H)$, we get that the normalization $\scrN_3$ of $\scrN$ in the ring of fractions of $\Sh(G_3,\scrX_3)/H_3$ together with the natural action of $G_3(\dbA_f^{(p)})$ on it, is the integral canonical model of $(G_3,\scrX_3,H_3)$. 

Let $W(\dbF)$ be the ring of Witt vectors with coefficients in an algebraic closure $\dbF$ of $\dbF_p$. We consider an arbitrary $\dbZ_{(p)}$-embedding $E(G_2,\scrX_2)_{(p)}\hookrightarrow W(\dbF)$. 

We can identify each connected component $\scrC_3$ of $\scrN_{3,W(\dbF)}$ with a connected component $\scrC$ of $\scrN_{W(\dbF)}$, cf. [Va7, Prop. 2.4.3 (c)]. Let $\scrC_2$ and $\scrC_1$ be the connected components of $\scrN_{2,W(\dbF)}$ and $\scrN_{1,W(\dbF)}$ (respectively) dominated by $\scrC_3$. The composite morphism $\scrC_3=\scrC\to\scrC_1$ of pro-finite covers, is a pro-\'etale cover (cf. Theorem 1.5 (iii)). Thus $\scrC_2$ is a pro-\'etale cover of $\scrC_1$. As the connected components of $\scrN_{2,W(\dbF)}$ are permuted transitively by $G_2(\dbA_f^{(p)})$ (cf. [Va1, Lem. 3.3.2]), by using $G_2(\dbA_f^{(p)})$-translates of $\scrC_2$ we get that $\scrN_{2,W(\dbF)}$ is a pro-\'etale cover of an open closed subscheme of $\scrN_{1,W(\dbF)}$. Thus $\scrN_2$ is a pro-\'etale cover of an open closed subscheme of $\scrN_1$. As $\scrN_1$ has the extension property, each closed subscheme of it which is flat over $E(G_1,\scrX_1)_{(p)}$ has also the extension property. From the last two sentences we get that the $E(G_2,\scrX_2)_{(p)}$-scheme $\scrN_2$ has the extension property, cf. [Va1, Rm. 3.2.3.1 6)]. 

It is easy to see that there exists a compact, open subgroup $\tilde H$ of $G_2(\dbA_f^{(p)})$ such that the morphism $\scrN_2\to\scrN_2/\tilde H$ is a pro-\'etale cover. As $\scrN$ is quasi-projective, we easily get that  $\scrN_2/\tilde H$ is a smooth, quasi-projective $E(G_2,\scrX_2)_{(p)}$-scheme. Thus $\scrN_2$ is the integral canonical model of $(G_2,\scrX_2,H_2)$ and it is quasi-projective. 

If $\scrN$ is projective, then $\scrN_1$ is projective (cf. Theorem 1.5 (iii)) and this implies that $\scrN_2$ is projective.\endproof

\bigskip
\noindent
{\boldsectionfont 4. The proof of the Main Theorem A}
\bigskip

In Subsection 4.1 we prove Theorem 1.6 (a) to (c). In Subsection 4.2 we prove Theorem 1.6 (d). Let $(G_1,\scrX_1)$ be a Shimura pair of abelian type. We assume that the group $G_{1,\dbQ_p}$ is unramified. Let $H_1$ be a hyperspecial subgroup of $G_1(\dbQ_p)$. Let $G_{1,\dbZ_p}$ be the reductive group scheme over $\dbZ_p$ such that its generic fibre is $G_{1,\dbQ_p}$ and we have $H_1=G_{1,\dbZ_p}(\dbZ_p)$. Let $H_1^{\ad}:=G_{1,\dbZ_p}^{\ad}(\dbZ_p)$; it is a hyperspecial subgroup of $G_1^{\ad}(\dbQ_p)$.

If the adjoint group $G_1^{\ad}$ is trivial, then Main Theorem A is well known (it is an easy consequence of either [Mi2, Rm. 2.16] or [Va1, Ex. 3.2.8]). Thus to prove the Main Theorem A we can assume that the adjoint group $G_1^{\ad}$ is non-trivial.

\bigskip\noindent
{\bf 4.1. Proofs of 1.6 (a) to (c).} Let $(G_1^{\ad},\scrX^{\ad},H_1^{\ad})=\prod_{i\in I} (G_i,\scrX_i,H_i)$ be the product decomposition into simple, adjoint Shimura triples. Thus each $(G_i,\scrX_i)$ is a simple, adjoint Shimura pair. Let $\scrN_i$ be the integral canonical model of $(G_i,\scrX_i,H_i)$, cf. Theorem 1.5 (ii). We consider the product $\scrN^{\ad}:=\prod_{i\in I} \scrN_{i,E(G_1^{\ad},\scrX_1^{\ad})_{(p)}}$ of $E(G_1^{\ad},\scrX_1^{\ad})_{(p)}$-schemes; it is the integral canonical model of $(G_1^{\ad},\scrX_1^{\ad},H_1^{\ad})$. Let $\scrN_1$ be the normalization of $\scrN^{\ad}$ in the ring of fractions of $\Sh(G_1,\scrX_1)/H_1$. We check that:

\medskip
{\bf (*)} the natural morphism $m_1:\scrN_1\to\scrN^{\ad}$ of $E(G_1^{\ad},\scrX_1^{\ad})_{(p)}$-schemes is pro-finite and a pro-\'etale cover of its image. 

\medskip
Let $f_3:(G_3,\scrX_3,H_3)\to (G_1,\scrX_1,H_1)$ be a cover such that at the level of reflex fields we have $E(G_3,\scrX_3)=E(G,\scrX)$ and the semisimple group cover $G_3^{\der}$ of $G_1^{\der}$ is the maximal one allowed by the abelian type, cf. [Va1, Rm. 3.2.7 10)]. Similarly we consider a cover $f_{3,i}:(G_{3,i},\scrX_{3,i},H_{3,i})\to (G_i,\scrX_i,H_i)$ such that at the level of reflex fields we have $E(G_{3,i},\scrX_{3,i})=E(G_i,\scrX_i)$ and the semisimple group cover $G_{3,i}^{\der}$ of $G_i$ is the maximal one allowed by the abelian type. The morphisms $\Sh(G_3,\scrX_3)/H_3\to\Sh(G_1,\scrX_1)/H_1$ and $\Sh(G_{3,i},\scrX_{3,i})/H_{3,i}\to\Sh(G_i,\scrX_i)/H_i$ are pro-\'etale covers, cf. [Va7,  Lem. 2.4.1]. In particular, we get that to check that the property (*) holds, we can assume that $G_1^{\der}$ is the maximal semisimple group cover of $G_1^{\ad}$ allowed by the abelian type. Let $(G_4,\scrX_4,H_4):=\prod_{i\in I} (G_{3,i},\scrX_{3,i},H_{3,i})$. We have $(G_4^{\ad},\scrX_4^{\ad},H_4^{\ad})=(G_1^{\ad},\scrX_1^{\ad},H_1^{\ad})$ and $G^{\der}_4=G_1^{\der}$. Based on [Va7, Prop. 2.4.3 (a) and (c)], to prove that the property (*) holds we can also assume that we have $(G_4,\scrX_4,H_4)=(G_1,\scrX_1,H_1)$. Thus to check that the property (*) holds, we can assume that the set $I$ has one element (i.e., $G_1^{\ad}$ is a simple, adjoint group over $\dbQ$). But this case follows from  Corollary 3.9. 

 As in the end of the proof of Corollary 3.9 we argue that $\scrN_1$ is the integral canonical model of $(G_1,\scrX_1,H_1)$ and it is quasi-projective. Thus Theorem 1.6 (a) holds. Theorem 1.6 (b) follows from the property (*) applied to the morphisms $m_1:\scrN_1\to\scrN^{\ad}$ and $m_2:\scrN_2\to\scrN^{\ad}$ of $E(G_1^{\ad},\scrX_1^{\ad})_{(p)}$-schemes, once we remark that $m_1$ is the composite of the functorial morphism $\scrN_2\to\scrN_1$ of $E(G_2,\scrX_2)_{(p)}$-schemes with $m_2$. 

Theorem 1.6 (c) follows from Theorem 1.6 (a) and Proposition 2.4.\endproof

\bigskip\noindent
{\bf 4.2. Proof of 1.6 (d).}  Let the following notations $(G_1^{\ad},\scrX^{\ad},H_1^{\ad})=\prod_{i\in I} (G_i,\scrX_i,H_i)$, $\scrN_i$, $\scrN^{\ad}:=\prod_{i\in I} \scrN_{i,E(G_1^{\ad},\scrX_1^{\ad})_{(p)}}$, and $\scrN_1$ be as Subsection 4.1. To prove that $\scrN_1$ is projective, it suffices to show that each $\scrN_i$ is projective. Thus we can assume that $G_1$ is a simple, adjoint group over $\dbQ$. Therefore we can appeal to the (notations of) Theorem 1.5. The connected components of $\scrN_1$ are permuted transitively by $G_1(\dbA_f^{(p)})$, cf. [Va1, Lem. 3.3.2]. Based on this and Theorem 1.5 (iii), to prove that $\scrN_1$ is projective it suffices to show that $\scrN$ is projective. Let $N\in\dbN\setminus (p\dbN\cup\{1,2\})$. Let $K(N)_p$ be the open closed subgroup of $\pmb{\text{GSp}}(L,\psi)(\dbA_f^{(p)})$ such that $K(N):=K_p\times K(N)_p$ is the subgroup of $\pmb{\text{GSp}}(L,\psi)(\widehat{\dbZ})$ formed by elements congruent to the identity modulo $N\widehat{\dbZ}$. Let $H(N)_p:=G(\dbA_f^{(p)})\cap K(N)_p$. Then $H(N):=H\times H(N)_p=G(\dbA_f)\cap K(N)$. 

As $N\ge 3$, a principally polarized abelian scheme with level-$N$ structure has no automorphism (see [Mu1, Ch. IV, 21, Thm. 5] for this result of Serre). This implies that $K(N)$ acts freely on $\scrA_{d,1,N}$. From this we get that $H(N)$ acts freely on $\scrN$. Therefore $\scrN$ is a pro-\'etale cover of the $E(G,\scrX)_{(p)}$-scheme $\scrN_N:=\scrN/H(N)$. The scheme $\scrN_N$ is the normalization of $\scrA_{d,1,N,\dbZ_{(p)}}$ in $\Sh(G,\scrX)/H(N)$ and therefore it is a finite $\scrA_{d,1,N,\dbZ_{(p)}}$-scheme. Therefore $\scrN_N$ is a quasi-projective $E(G,\scrX)_{(p)}$-scheme. Thus the $E(G,\scrX)_{(p)}$-scheme $\scrN_N$ is projective if and only if $\scrN$ is projective.

We write $G_1=\Res_{F/\dbQ} J_1$, where $F_1$ is a totally real number field and $J_1$ is an absolutely simple adjoint group over $F_1$ (cf. Subsection 2.2). If $(G_1,\scrX_1)$ has compact factors, then the projectiveness of $\scrN_N$ is implied by [Va6, Cor. 4.3] and [Va15, Part I, Lem. 2.2.4]. Thus we can assume that each simple factor of $G_{1,\dbR}$ is non-compact. From this and the fact that $G_1$ has $\dbQ$--rank $0$, one gets that there exists a finite prime $w$ of $F$ such that $J_{1,F_w}$ is anisotropic i.e., has $F_w$-rank $0$ (cf. [Lee, Thm. 2.5]; here $F_w$ is the completion of $F$ at $w$). Thus, as $G_1$ has $\dbQ$--rank $0$, from the classification of absolutely simple classical adjoint groups over number fields (see [Ti1]) and of simple adjoint Shimura pairs of abelian type (see Subsection 1.9) we get that $(G_1,\scrX_1)$ if unitary (cf. [Va6, Rm. 2.3.2]). The fact that $\scrN_N$ is projective is equivalent to the following property (to be compared with [Va6, Prop. 2.7 (b) and Cor. 4.3] and [Va15, Part I, Lem. 2.2.4]):

\medskip\noindent
$(\flat)$ For each point of $\scrN_N$ with values in a number field $E$ that contains $E(G,\scrX)$, the abelian variety $\scrD$ over $E$ associated naturally to the point has potentially good reduction with respect to all primes of $E$ dividing $p$ provided the adjoint of the Mumford--Tate group $H_{\scrD}$ of a fixed pull-back of $\scrD$ to $\dbC$ is $G_1$ itself. 

\medskip
There are at least four ways to prove that $\scrN_N$ is projective (equivalently that the property $(\flat)$ holds) to be listed below using bullets.

\medskip
$\bullet$ As $(G_1,\scrX_1)$ if unitary, to check that $\scrN_N$ is projective (equivalently that the property $(\flat)$ holds) we can assume that we are in the Case 1 of Subsection 3.1 (cf. [Va8, Cor. 4.10]) and in particular that the injective map $(G,\scrX)\hookrightarrow (GSp(W,\psi),\scrS)$ is a PEL type embedding. Thus $\scrD$ has potentially good reduction everywhere if $H_{\scrD}^{\ad}=G_1$ (cf. [Pa, Prop. 4.2.13]) and therefore the property $(\flat)$ holds. 

\smallskip
$\bullet$ As in the first bullet, we can assume that the injective map $(G,\scrX)\hookrightarrow (GSp(W,\psi),\scrS)$ is a PEL type embedding. Thus the fact that $\scrN_N$ is projective follows as well from adapting the arguments of [Va6, Subsect. 3.4] to the context in which a non-archimedean (instead of an archimedean) completion $F_w$ of $F$ is such that $j_{1,F_w}$ is anisotropic.   

\smallskip
$\bullet$ As in the first bullet, we can assume that the injective map $(G,\scrX)\hookrightarrow (GSp(W,\psi),\scrS)$ is a PEL type embedding. Thus the fact that $\scrN_N$ is projective follows as well from Mumford's constructions of [Mu2], cf. [Lan].   

$\bullet$ The fact that the property $(\flat)$ holds is checked in [Lee] using a refinement of [Pa, Thm. 1.6.1] (which in essence is only a variant of [Pa, Prop. 4.2.13]).

\medskip
Based on either one of these four bullets one concludes that $\scrN_N$ is projective. Thus Theorem 1.6 (d) holds. This ends the proof of the Main Theorem A.\endproof

\bigskip
\noindent
{\boldsectionfont 5. The proof of the Main Theorem B}
\bigskip

In this section we prove the Main Theorem B. The fact that the condition 1.7 (i) holds is implied by Theorem 1.6 (a) and Proposition 2.4. Thus we have $\scrN=\scrN^{\s}$, cf. notations of Subsubsection 2.2.1. Let $(v_{\alpha})_{\alpha\in\scrJ}$, $k$, $W(k)$, $B(k)$, and $\sigma$ be as in Subsubsection 2.2.1. Let $k(v)$ be the residue field of $v$. Let $O:=O_{(v)}$. We have $\scrN^{\s}_{O}(W(k))=\scrN_{O}(W(k))$, cf. [Va15, Thm. 1.5 (a)]. 

For $z\in\scrN^{\s}_{O}(W(k))$, let $(A,\lambda_A):=z^*(\scrA,\Lambda_{\scrA})$. Let $(M,F^1,\phi,\psi_M)$ be the principally quasi-polarized filtered Dieudonn\'e module of the principally quasi-polarized  $p$-divisible group $(D,\lambda_D)$ of $(A,\lambda_A)$. For each $\alpha\in\scrJ$ let $t_{\alpha}$ (resp. $u_{\alpha}$) be the de Rham component (resp. the $p$-component of the \'etale component) of the Hodge cycle on $A_{B(k)}$ that corresponds naturally to $v_{\alpha}$. Each $t_{\alpha}$ belongs to $\scrT(M)[{1\over p}]$ and it is fixed by $\phi$ (cf. [Va10, Thm. 5.1.6 and Cor. 5.1.7]). We refer to $(M,F^1,\phi,(t_{\alpha})_{\alpha\in\scrJ})$ as the {\it filtered Dieudonn\'e module with tensors} attached to $z$ (to be compared with [Va10, Subsubsect. 5.1.8]). Let $\scrG$ be the schematic closure in $\pmb{\text{GL}}_M$ of the subgroup of $\pmb{\text{GL}}_{M[{1\over p}]}$ that fixes $t_{\alpha}$ for all $\alpha\in\scrJ$. If $\scrG$ is a reductive group scheme over $W(k)$, then we call $(M,F^1,\phi,\scrG)$ as the {\it filtered Shimura $F$-crystal} attached to $z$ (see [Va5,9,10,14] for the notion filtered Shimura $F$-crystal). We recall from [Va15, Subsubsect. 3.5.1 and Lem. 2.3.4 (a)] that $\scrN_O^{\m}$ is the $G(\dbA_f^{(p)})$-invariant, open subscheme of $\scrN^{\s}_O$ with the property that a point $z\in\scrN^{\s}_{O}(W(k))$ factors through $\scrN^{\m}_{O}$ if and only if there exists an isomorphism
$$\rho_D:(M,(t_{\alpha})_{\alpha\in\scrJ})\arrowsim (L_p^{\vee}\otimes_{\dbZ_p} W(k),(v_{\alpha})_{\alpha\in\scrJ}).$$
Thus the condition 1.7 (ii) holds if and only if the isomorphism $\rho_D$ exists for each  $z\in\scrN^{\s}_{O}(W(k))$. Let $\mu:\dbG_m\to \scrG$ be a cocharacter as in [Va15, Subsect. 3.2]; let $M=F^1\oplus F^0$  be a direct sum decomposition such that $\dbG_m$ acts through $\mu$ on each $F^i$ via the weight $-i$. 

Based on [Va15, Thm. 3.2.2 (a)], it suffices to prove the existence of the isomorphism $\rho_D$ in the case when $p=2$ and $D$ is not a direct sum of \'etale and connected $2$-divisible groups. This implies that the $2$-rank of $A_k$ is positive. Based on [Va15, Lem. 2.3.4 (a)], the $2$-adic Galois representation associated to $A_{B(k)}$ can be identified with a homomorphism
$$\varrho_D:\Gal(B(k))\to G_{\dbZ_2}(\dbZ_2).$$
We consider the induced Galois representation
 $$\varrho_D^{\ad}:\Gal(B(k))\to G_{\dbZ_2}^{\ad}(\dbZ_2).$$
\indent
We consider a connected component $\scrC$ of $\scrN_{O}$ and a connected component $\bar\scrC$ of $\scrC_{k(v)}$. As the connected components of $\scrN_{O}$ are permuted transitively by $G(\dbA_f^{(2)})$ (cf. [Va1, Lem. 3.3.2]) and $\scrN^{\m}_{O}$ is a $G(\dbA_f^{(2)})$-invariant, open subscheme of $\scrN_{O}$, to prove that the isomorphism $\rho_D$ exists we can assume that $z\in\scrC(W(k))$. As the special fibre $\scrN^{\m}_{k(v)}$ of $\scrN^{\m}_{O}$ is an open closed subscheme of the special fibre $\scrN_{k(v)}$ of $\scrN_{O}$ (cf. [Va15, Thm. 1.7 (a)]), it suffices to show that there exists {\it one} point $z\in \scrC(W(k))$ that lifts a $k$-valued point of $\bar\scrC$ for which the isomorphism $\rho_D$ exists. Based on [Va15, Thm. 3.2.2 (a)] it suffices to show that there exists {\it one} point $z\in \scrC(W(k))$ that lifts a $k$-valued point of $\bar\scrC$ for which $\im(\varrho_D)$  is contained in the group of $\dbZ_2$-valued points of a torus of $G_{\dbZ_2}$. Thus it suffices to show:

\medskip
{\bf ($\sharp$)} given a connected component $\bar\scrC$ of the special fibre of $\scrN_{O}$, there exists {\it one} point $z\in \scrN_{O}(W(k))$ that lifts a $k$-valued point of $\bar\scrC$ for which $\im(\varrho_D^{\ad})$  is contained in the group of $\dbZ_2$-valued points of a torus of $G_{\dbZ_2}$. 

\medskip
But ($\sharp$) is a condition that is invariant under the operations of taking Hodge twists and Hodge quasi products introduced in [Va8, Subsects. 2.4 and 5.3])  and thus (to be compared with the shifting process [Va8, Prop. 5.4.1]) based on Theorem 1.6 and its proof in Subsection 4.1, we can assume that $(G_1,\scrX_1):=(G^{\ad},\scrX^{\ad})$ is a simple adjoint Shimura pair and that we are in the context of Theorem 1.5. Thus we can assume that we are in one of the six cases of Subsections 3.1 to 3.6 and we are left to check that the condition ($\sharp$) holds. 

First, we will check case by case that the isomorphism $\rho_D$ always exists i.e., we have $\scrN^{\m}_{O}=\scrN_{O}$. In the first three cases the existence of the isomorphism $\rho_D$ is well known (for instance, it follows from [Ko, Lem. 7.2] and [Va11, Fact 2.5.1 (a)]). In the fourth case the existence of the isomorphism $\rho_D$ follows from [Va14, Cor. 5.5]. In the fifth (resp. sixth) case the identity $\scrN^{\m}_{O}=\scrN_{O}$ follows from [Va15, Thm. 1.7 (c)] (resp. from [Va15, Thm. 1.7 (c)] and Proposition 3.6.1). Thus in all six cases we have  $\scrN^{\m}_{O}=\scrN_{O}$.

Second, we will check that the identity $\scrN^{\m}_{O}=\scrN_{O}$ implies that the condition ($\sharp$) holds. As the isomorphism $\rho_D$ exists, $\scrG$ is a reductive group scheme over $W(k)$ isomorphic to $G_{\dbZ_{(p)}}\times_{\Spec\dbZ_{(p)}} \Spec W(k)$. From this and the identity $\scrN_{O}=\scrN^{\s}_{O}$ we get that the the triple $(f,L,v)$ is a standard Hodge situation in the sense of either [Va10, Def. 5.1.2] or [Va13, Def. 1.5.3]. To check that the condition ($\sharp$) holds we can assume that $k=\dbF$ and we consider a point $z\in \scrN_O(W(\dbF))$ that lifts an arbitrary point $y\in\bar\scrC(\dbF)$. We claim that we can choose $y\in\bar\scrC(\dbF)$ such that we can assume that there exist a maximal torus $\scrT$ of $\scrG$ and a Borel subgroup scheme $\scrB$ of $\scrG$ that contains $\scrT$, such that the following two properties hold:

\medskip
{\bf (i)} The cocharacter $\mu$ factors through $\scrT$.

\smallskip
{\bf (ii)} The Lie algebra $\Lie(\scrT)$ (resp. $\Lie(\scrB)$) is normalized (resp. left invariant) by $\phi$.

\medskip
This claim is a direct consequence of [Va13, Thm. 1.6.1, Cor. 3.1.9, and Subsect. 4.1.1] (more precisely, it suffices to choose $y\in\scrC(\dbF)$ to be a Shimura-ordinary point in the sense of [Va13, Def. 1.6.4]). 

Here is a second way to argue that the claim holds.
Based on [Va15, property 3.5.1 (iii)], the filtered Dieudonn\'e module with tensors attached to another point $z_1\in\scrN_O(W(\dbF))$ that lifts a point $y_1\in\bar\scrC(\dbF)$ is isomorphic to $(M,F^1,g_1\phi,(t_{\alpha})_{\alpha\in\scrJ})$ for some element $g_1\in \scrG(W(\dbF))$. Thus from [Va9, Basic Thm. 12.2 and Rm. 12.4 (a)] we get the existence of a level $1$ stratification of $\bar\scrC$ with the property that two points $y_1,y_2\in\bar\scrC(\dbF)$ belong to the same stratum if and only if, up to isomorphisms, we can assume that the elements $g_1$ and $g_2$ that correspond naturally to $y_1$ and $y_2$ (respectively) are congruent modulo $p$. This level $1$ stratification has an open, Zariski dense stratum and to ease the notations we can assume that $y$ belongs to it. Thus from [Va9, Ex. 5.6 and Basic Thm. D] we get that there exists an element $h\in\Ker(\scrG(W(\dbF))\to\scrG(\dbF))$, a maximal torus $\scrT$ of $\scrG$, and a Borel subgroup scheme $\scrB$ of $\scrG$ that contains $\scrT$, such that the property (i) holds and moreover the Lie algebra $\Lie(\scrT)$ (resp. $\Lie(\scrB)$) is normalized (resp. left invariant) by $h\phi$. Based on [Va5, Prop. 4.3.1] and [Va9, Ex. 5.6], we get that there exists an element $h_1\in\scrG(W(\dbF))$ such that we have $h_1h\phi =\phi h_1$  i.e., up to an isomorphism and up to a replacement of $z$ by another point $z_{h_1}\in\scrN_O(W(\dbF))$ that lifts $y$ we can assume that $h=1_M$ (it is [Va15, Lem. 3.5.2] that guarantees that we can choose $z_{h_1}$ such that the Hodge filtration of $M$ defined by $z_{h_1}^*(\scrA)$ is $h_1(F^1)$; note that $h_1$ modulo $2$ normalizes $F^1/2F^1$). As $h=1_M$, we get that the property (ii) holds as well and this ends the second way to argue that the claim holds.

From the property (i) we get that  $(M,F^1,\phi)$ is a direct sum of filtered Dieudonn\'e module over $k$ whose associated Dieudonn\'e modules are isoclinic (to be compared with [Va11, Susubsect. 4.1.2]). Let $\scrT_{\dbQ_2}^{\acute{et}}$ be the torus of $G_{\dbQ_2}$ which corresponds to the torus $\scrT_{B(\dbF)}$ of $\scrG_{B(\dbF)}$ via the Fontaine comparison theory for $D$, cf. [Va15, Appendix, Subsect. B5] and properties (i) and (ii). Let $\scrT^{\acute{et}}$ be the schematic closure of $\scrT_{\dbQ_2}^{\acute{et}}$ in $G_{\dbZ_2}$; we have $\im(\varrho_D)\leqslant \scrT^{\acute{et}}(\dbZ_2)$. 

Thus the condition ($\sharp$) holds provided we can modify the lift $z$ of $y$ in such a way that $F^1$ does not change and moreover $\scrT^{\acute{et}}$ is a torus over $\dbZ_2$. From [Va15, Lem. 3.5.2] and the proof of [Va15, Appendix, Thm. B7] we get that we can choose the lift $z$ of $y$ such that $F^1$ does not change and $D$ is a direct sum of $2$-divisible groups whose special fibres are isoclinic. Based on this, either as in [Va11, Lem. 4.1.1] we argue that $\scrT^{\acute{et}}$ is a torus or directly from [Va11, Cor 1.3 (a)] applied to $D$ and $(M,F^1,\scrT,(t_{\alpha})_{\alpha\in\scrJ_{\scrT}})$ we get that $\scrT^{\acute{et}}$ is a torus; here $(t_{\alpha})_{\alpha\in\scrJ_{\scrT}}$ is an arbitrary family of tensors of $\scrT(M[{1\over 2}])$ fixed by $\phi$ and $\scrT_{B(k)}$ with the property that $\scrT_{B(k)}$ is the subgroup of $\pmb{GL}_{M[{1\over 2}]}$ that fixes $t_{\alpha}$ for all $\alpha\in\scrJ_{\scrT}$ (it exists, cf. [Va11, Lem. 2.5.3] and properties (i) and (ii)). 

Thus ($\sharp$) and the property 1.7 (ii) hold. This ends the proof of the Main Theorem B.\endproof

\medskip\noindent
{\bf Acknowledgments.} We would like to thank Binghamton University and IAS, Princeton for good working conditions. This research was partially supported by the NSF grants  DMF 97-05376 and DMS \#0900967. 

\bigskip
\references{37}
{\nspace{

\Ref[Ar1] M. Artin, 
\sl Algebraization of formal moduli. I, 
\rm Global Analysis (Papers in Honor of K. Kodaira), 21--71, Univ. Tokyo Press, Tokyo, 1969.

\Ref[Ar2] M. Artin, 
\sl Versal deformations and algebraic stacks, 
\rm Invent. Math. {\bf 27} (1974), 165--189.

\Ref[BB]
W. Baily and A. Borel,
\sl Compactification of arithmetic quotients of bounded
symmetric domains,
\rm Ann. of Math. (2) {\bf 84} (1966), no. 3, 442--528.

\Ref[BLR]
S. Bosch, W. L\"utkebohmert, and M. Raynaud,
\sl N\'eron models,
\rm Ergebnisse der Mathematik und ihrer Grenzgebiete (3), Vol. {\bf 21}, Springer-Verlag, Berlin, 1990.

\Ref[De1]
P. Deligne,
\sl Travaux de Shimura,
\rm S\'eminaire  Bourbaki, 23\`eme ann\'ee (1970/71), Exp. No. 389, Lecture Notes in Math., Vol. {\bf 244}, 123--165, Springer-Verlag, Berlin, 1971.

\Ref[De2]
P. Deligne,
\sl Vari\'et\'es de Shimura: interpr\'etation modulaire, et
techniques de construction de mod\`eles canoniques,
\rm Automorphic forms, representations and $L$-functions (Oregon State Univ., Corvallis, OR, 1977), Part 2,  247--289, Proc. Sympos. Pure Math., {\bf 33}, Amer. Math. Soc., Providence, RI, 1979.

\Ref[De3]
P. Deligne,
\sl Hodge cycles on abelian varieties,
\rm Hodge cycles, motives, and Shimura varieties, Lecture Notes in Math., Vol.  {\bf 900}, 9--100, Springer-Verlag, Berlin-New York, 1982.

\Ref[Fa]
G. Faltings,
\sl Integral crystalline cohomology over very ramified
valuation rings,
\rm J. Amer. Math. Soc. {\bf 12} (1999), no. 1, 117--144.

\Ref[FC]
G. Faltings and C.-L. Chai, 
\sl Degeneration of abelian varieties, 
\rm Ergebnisse der Mathematik und ihrer Grenzgebiete (3), Vol. {\bf 22}, Springer-Verlag, Berlin, 1990.

\Ref[Ha] 
R. Hartshorne, 
\sl Algebraic geometry, 
\rm Grad. Texts in Math., {\bf 52}, Springer-Verlag, Berlin, 1977.

\Ref[Ki1] M. Kisin,
\sl Crystalline representations and $F$-crystals, 
\rm Algebraic geometry and number theory,  459--496, Progr. Math., Vol. {\bf 253}, Birkh\"auser,
  Boston, MA, 2006.

\Ref[Ki2] M. Kisin,
\sl Modularity of 2-adic Barsotti--Tate representations, \rm Invent. Math.. Vol. {\bf 178} (2009), no. 3,  587--634.

\Ref[Ki3] 
M. Kisin,
\sl Integral canonical models of Shimura varieties of abelian type, 
\rm J. Amer. Math. Soc. {\bf 23} (2010), no. 4,  967--1012.

\Ref[Ko]
R. E. Kottwitz,
\sl Points on some Shimura varieties over finite fields,
\rm J. Amer. Math. Soc. {\bf 5} (1992), no. 2, 373--444.

\Ref[Lan]
K.-W. Lan,
\sl Elevators for degenerations of PEL structures,
\rm Math. Res. Lett. {\bf 18} (2011), 889--907.

\Ref[Lee]
D. U. Lee,
\sl A proof of a conjecture of Morita,
\rm 10 pages, to appear in Bull. London Math. Soc., 

http://blms.oxfordjournals.org/content/early/2012/06/28/blms.bdr104.short?rss=1. 

\Ref[LR]
R. Langlands and M. Rapoport,
\sl Shimuravariet\"aten und Gerben, 
\rm J. Reine Angew. Math. {\bf 378} (1987), 113--220.

\Ref[Mi1]
J. S. Milne,
\sl Canonical models of (mixed) Shimura varieties and automorphic vector bundles, 
\rm  Automorphic Forms, Shimura varieties and L-functions, Vol. I (Ann Arbor, MI, 1988), 283--414, Perspectives in Math., Vol. {\bf 10}, Academic Press, Inc., Boston, MA, 1990.

\Ref[Mi2]
J. S. Milne,
\sl The points on a Shimura variety modulo a prime of good
reduction,
\rm The Zeta functions of Picard modular surfaces, 153--255, Univ. Montr\'eal, Montreal, Quebec, 1992.

\Ref[Mi3]
J. S. Milne,
\sl Shimura varieties and motives,
\rm Motives (Seattle, WA, 1991), Part 2, 447--523, Proc. Sympos. Pure Math., Vol. {\bf 55}, Amer. Math. Soc., Providence, RI, 1994.

\Ref[Mi4]
J. S. Milne,
\sl Descent for Shimura varieties,
\rm Michigan Math. J. {\bf 46} (1999), no. 1, 203--208.

\Ref[Mi5]
J. S. Milne,
\sl Points on Shimura varieties over finite fields: the conjecture of Langlands and Rapoport,

\rm 40 pages manuscript, November 11, 2009, http://arxiv.org/abs/0707.3173.

\Ref[Mo]
Y. Morita,
\sl On potential good reduction of abelian varieties,
\rm J. Fac. Sci. Univ. Tokyo Sect. I A Math. {\bf 22} (1975), no. 3, 437--447.

\Ref[Moo]
B. Moonen,
\sl Models of Shimura varieties in mixed characteristics,
\rm Galois representations in arithmetic algebraic geometry (Durham, 1996),  267--350, London Math. Soc. Lecture Note Ser., {\bf 254}, Cambridge Univ. Press, Cambridge, 1998.  

\Ref[Mu1]
D. Mumford,
\sl Abelian varieties,
\rm Tata Inst. of Fund. Research Studies in Math., No. {\bf 5}, Published for the Tata Institute of Fundamental Research, Bombay; Oxford Univ. Press, London, 1970, reprinted 1988.

\Ref[Mu2]
D. Mumford,
\sl An analytic construction of degenerating abelian varieties over complete rings,
\rm Compos. Math. {\bf 24} (1972), no. 3, 239--272.

\Ref[MFK]
D. Mumford, J. Fogarty, and F. Kirwan,
\sl Geometric invariant theory. Third enlarged edition,
\rm Ergebnisse der Mathematik und ihrer Grenzgebiete (2), Vol. {\bf 34}, Springer-Verlag, Berlin, 1994. 

\Ref[MS]
J. S. Milne and K.-Y. Shih
\sl Conjugates of Shimura varieties,
\rm  Hodge cycles, motives, and Shimura varieties, Lecture Notes in Math., Vol.  {\bf 900}, 280--356, Springer-Verlag, Berlin-New York, 1982.

\Ref[No] 
R. Noot, 
\sl Models of Shimura varieties in mixed characteristic, 
\rm J. Algebraic Geom. {\bf 5} (1996), no. 1, 187--207.

\Ref[Pa]
F. Paugam,
\sl Galois representations, Mumford--Tate groups and good reduction of abelian varieties,
\rm Math. Ann. {\bf 329} (2004), no. 1, 119--160. Erratum: Math. Ann. {\bf 332} (2004), no. 4, 937.

\Ref[PY]
G. Prasad and J.-K. Yu,
\sl On quasi-reductive group schemes. With an appendix by Brian Conrad,
\rm J. Algebraic Geom. {\bf 15} (2006), no. 3,  507--549. 

\Ref[Sa1]
I. Satake,
\sl Holomorphic imbeddings of symmetric domains into a Siegel space,
\rm  Amer. J. Math. {\bf 87} (1965), 425--461.

\Ref[Sa2]
I. Satake,
\sl Symplectic representations of algebraic
groups satisfying a certain analyticity condition,
\rm Acta Math. {\bf 117} (1967), 215--279.

\Ref[Ti1]
J. Tits, 
\sl Classification of algebraic semisimple groups, 
\rm Algebraic Groups and Discontinuous Subgroups (Boulder, CO, 1965),  33--62, Proc. Sympos. Pure Math., Vol. {\bf 9}, Amer. Math. Soc., Providence, RI, 1966.

\Ref[Ti2]
J. Tits,
\sl Reductive groups over local fields, 
\rm Automorphic forms, representations and $L$-functions (Oregon State Univ., Corvallis, OR, 1977), Part 1,  29--69, Proc. Sympos. Pure Math., Vol. {\bf 33}, Amer. Math. Soc., Providence, RI, 1979.

\Ref[Va1]
A. Vasiu,
\sl Integral canonical models of Shimura varieties of preabelian type,
\rm Asian J. Math. {\bf 3} (1999), no. 2, 401--518.

\Ref[Va2] 
A. Vasiu, 
\sl Surjectivity criteria for $p$-adic representations, Part I,
\rm Manuscripta Math. {\bf 112} (2003), no. 3, 325--355.

\Ref[Va3] 
A. Vasiu,
\sl A purity theorem for abelian schemes,
\rm Michigan Math. J. {\bf 54} (2004), no. 1, 71--81.

\Ref[Va4] 
A. Vasiu,
\sl  On two theorems for flat, affine groups schemes over a discrete valuation ring,
\rm Centr. Eur. J. Math. {\bf 3} (2005), no. 1, 14--25.

\Ref[Va5]
A. Vasiu,
\sl Crystalline boundedness principle,
\rm Ann. Sci. \'Ecole Norm. Sup. {\bf 39} (2006), no. 2, 245--300.

\Ref[Va6] 
A. Vasiu,
\sl Projective integral models of Shimura varieties of Hodge type with compact factors,
\rm J. Reine Angew. Math. {\bf 618} (2008), 51--75.

\Ref[Va7] 
A. Vasiu,
\sl Integral canonical models of unitary Shimura varieties,
\rm Asian J. Math. {\bf 12} (2008), no. 2,  151--176.

\Ref[Va8]
A. Vasiu,
\sl Some cases of the Mumford-Tate conjecture and Shimura varieties,   
\rm Indiana Univ. Math. J. {\bf 57} (2008), no. 1, 1--76.

\Ref[Va9] 
A. Vasiu,
\sl Mod $p$ classification of Shimura $F$-crystals,
\rm Math. Nachr. {\bf 283} (2010), no. 8, 1068--1113.

\Ref[Va10]
A. Vasiu,
\sl Manin problems for Shimura varieties of Hodge type,
\rm J. Ramanujan Math. Soc. {\bf 26} (2011), no. 1, 31--84.

\Ref[Va11] 
A. Vasiu,
\sl A motivic conjecture of Milne,
\rm 67 pages to appear in J. Reine Angew. Math., http://www.degruyter.com/view/j/crelle.ahead-of-print/crelle-2012-0009/crelle-2012-0009.xml?format=INT.

\Ref[Va12] 
A. Vasiu,
\sl Integral models in unramified mixed characteristic $(0,2)$ of hermitian orthogonal Shimura varieties of PEL type, Part I,
\rm 44 pages, to appear in J. Ramanujan Math. Soc., available at http://arxiv.org/abs/math/0307205.

\Ref[Va13] 
A. Vasiu,
\sl Generalized Serre--Tate ordinary theory,
\rm 196 pages (including contents and index), to be published by International Press, Inc., all copyrights reserved to International Press, Inc., http://www.math.binghamton.edu/adrian/\#reductive.

\Ref[Va14] 
A. Vasiu,
\sl Integral models in unramified mixed characteristic $(0,2)$ of hermitian orthogonal Shimura varieties of PEL type, Part II,
\rm 24 pages manuscript, June 18, 2012, http://www.math.binghamton.edu/adrian/.

\Ref[Va15]
A. Vasiu,
\sl Good Reductions of Shimura varieties of hodge type in arbitrary unramified mixed characteristic, Part I,

\rm 51 pages manuscript, July 24, 2012, http://xxx.arxiv.org/abs/0707.1668.

\Ref[Va16]
A. Vasiu,
\sl On the Tate and Langlands--Rapoport conjectures for special fibres of integral canonical model of Shimura varieties of abelian type,
\rm 51 pages manuscript, July 24, 2012, http://www.math.binghamton.edu/adrian/.

\Ref[VZ]
A. Vasiu and T. Zink,
\sl Purity results for $p$-divisible groups and abelian schemes over regular bases of mixed characteristic, 
\rm Documenta Math. {\bf 15} (2010), 571--599.

\Ref[Zi]
T. Zink,
\sl Isogenieklassen von Punkten von Shimuramannigfaltigkeiten mit Werten in einem endlichen K\"orper,
\rm Math. Nachr. {\bf 112} (1983), 103--124.

}}
\medskip
\hbox{Adrian Vasiu,}
\hbox{Department of Mathematical Sciences, Binghamton University,}
\hbox{P. O. Box, Binghamton, New York 13902-6000, U. S. A.}
\hbox{e-mail: adrian\@math.binghamton.edu,\;\;fax: 1-607-777-2450,\;\;phone 1-607-777-6036}

\enddocument